\documentclass[oneside, reqno]{amsart}
\pdfoutput=1

\usepackage{subfiles}
\usepackage[utf8]{inputenc}
\usepackage[dvipsnames]{xcolor}
\usepackage{amssymb}
\usepackage{amsmath}
\usepackage{amsthm}
\usepackage[colorinlistoftodos]{todonotes}
\usepackage{microtype}
\usepackage{mathtools}
\usepackage{tikz}
\usepackage[bookmarks, colorlinks, allcolors=Blue]{hyperref}
\usepackage{tikz-cd}
\usepackage{lipsum}
\usepackage[foot]{amsaddr}

\usepackage[shortlabels]{enumitem}
\setlist[itemize]{leftmargin=2\parindent}
\setlist[enumerate]{label=(\arabic*)}

\usepackage[capitalize, noabbrev, nameinlink]{cleveref}
\Crefname{equation}{}{}
\Crefname{enumi}{}{}
\Crefname{axiom}{Axiom}{Axioms}

\usepackage[
    backend=biber,
    style=numeric-comp,
    giveninits=true,
    maxbibnames=10,
    minbibnames=10,
    mincrossrefs=10,
    bibencoding=utf8]{biblatex} % backref

\AtBeginBibliography{
    \raggedright
    \footnotesize
}

\ifSubfilesClassLoaded{%
    \AtEndDocument{%
        \printbibliography
    }
}{}

\makeatletter
\renewcommand{\@mkboth}[2]{}
\makeatother
\usetikzlibrary{decorations.markings}
\usetikzlibrary{cd}

\tikzcdset{
    diagrams={ampersand replacement=\&},
    postmark/.style={
        postaction={
            decorate,
            decoration={
                markings,
                mark={#1}
            }
        }
    },
    dagger arrow/.style={
        postmark=at position 0.5 with {
            \draw[-] (#1,1.5pt) -- (#1,-1.5pt);
        },
        labels={
            inner sep=0.75ex,
        }
    },
    mono/.style={tail},
    epi/.style={two heads},
    dagger mono/.style={dagger arrow={0pt}, mono},
    dagger epi/.style={dagger arrow={-1pt}, epi},
    left adjoint/.style={
        phantom,
        "\dashv"{font=\scriptsize},
        sloped,
        rotate={90}
    }
}
\input{preamble/widebar}
\theoremstyle{plain}
\newtheorem{theorem}{Theorem}
\newtheorem*{theorem*}{Theorem}

\newtheorem{lemma}[theorem]{Lemma}
\newtheorem{proposition}[theorem]{Proposition}
\newtheorem{corollary}[theorem]{Corollary}

\theoremstyle{definition}
\newtheorem{example}[theorem]{Example}
\newtheorem{definition}[theorem]{Definition}
\newtheorem{axiom}{Axiom}
\makeatletter
\newcommand{\neutralize}[1]{\expandafter\let\csname c@#1\endcsname\count@}
\makeatother
\newenvironment{alternativeaxiom}[1]
  {\renewcommand{\theaxiom}{\ref*{#1}$'$}%
   \neutralize{axiom}\phantomsection
   \begin{axiom}}
  {\end{axiom}}

\theoremstyle{remark}
\newtheorem{remark}[theorem]{Remark}
\newtheorem*{remark*}{Remark}

\renewcommand{\leq}{\leqslant}
\renewcommand{\geq}{\geqslant}
\newcommand{\cleq}{\preccurlyeq}

\newcommand{\nqquad}{\mkern-36mu}

\newcommand{\blank}{\underline{\phantom{m}}}

\DeclarePairedDelimiter{\paren}{(}{)}
\DeclarePairedDelimiter{\norm}{\|}{\|}
\DeclarePairedDelimiter{\abs}{|}{|}
\DeclarePairedDelimiterX{\innerProd}[2]{\langle}{\rangle}{#1 \delimsize\vert #2}
\DeclarePairedDelimiterX{\setb}[2]{\lbrace}{\rbrace}{#1 \,\delimsize\vert\,\mathopen{} #2}
\DeclarePairedDelimiter{\set}{\lbrace}{\rbrace}

\newcommand{\Disk}{\mathcal{D}}

\newcommand{\NZDisk}{\Disk_{\times}}
\newcommand{\SAScalars}{\mathcal{R}}
\newcommand{\PosScalars}{\mathcal{R}_+}
\newcommand{\Scalars}{\mathcal{C}}

\newcommand{\field}[1]{\mathbb{#1}}
\newcommand{\Nats}{\field{N}}
\newcommand{\Rats}{\field{Q}}
\newcommand{\PosRats}{\Rats_+}
\newcommand{\TropReals}{\field{TR}_+}
\newcommand{\Reals}{\field{R}}
\newcommand{\PosReals}{\Reals_+}
\newcommand{\Comps}{\field{C}}

\newcommand{\cat}[1]{\mathbf{#1}}
\newcommand{\Con}{\cat{Con}}
\newcommand{\FCon}{\cat{FCon}}
\newcommand{\Hilb}{\cat{Hilb}}
\newcommand{\FHilb}{\cat{FHilb}}
\newcommand{\Set}{\cat{Set}}
\newcommand{\C}{\cat{C}}

\newcommand{\D}{\cat{D}}
\newcommand{\A}{\cat{A}}

\DeclareMathOperator{\Par}{Par}
\newcommand{\normSq}{\mathcal{N}}

\newcommand{\pair}[2]{\begin{bsmallmatrix} #1 \\ #2 \end{bsmallmatrix}}
\newcommand{\copair}[2]{\begin{bsmallmatrix} #1 & #2 \end{bsmallmatrix}}

\newcommand{\id}[1][]{\ensuremath{{1_{#1}}}}

\DeclareMathOperator{\Ker}{Ker}
\DeclareMathOperator{\Coker}{Coker}
\DeclareMathOperator{\coker}{coker}
\DeclareMathOperator{\colim}{colim}

\renewcommand{\Im}{\operatorname{Im}}

\DeclareMathOperator{\Coim}{Coim}

\newcommand{\inl}{i_1}
\newcommand{\inr}{i_2}

\newcommand{\shortArrow}[3][]{\begin{tikzcd}[cramped,sep=small]#2 \arrow[r,#1] \& #3 \end{tikzcd}}

\newcommand{\daggerMono}{\shortArrow[dagger mono]}
\newcommand{\daggerEpi}{\shortArrow[dagger epi]}
\newcommand{\mono}{\shortArrow[mono]}
\newcommand{\epi}{\shortArrow[epi]}

\newlength{\UpperCaseHeight}
\AtBeginDocument{\settoheight{\UpperCaseHeight}{A}}

\newcommand*{\tikzdown}[1][]{\mathbin{\begin{tikzpicture}[commutative diagrams/every diagram, baseline=0]
\coordinate (T) at (0,-0.4ex);
\coordinate (S) at (0,\UpperCaseHeight+0.2ex);
\draw[commutative diagrams/.cd,
        every arrow,
        every label,
        #1]
    (S) -- (T);
\end{tikzpicture}}}

\newcommand*{\CommaCat}[3][]{#2 \tikzdown[#1] #3}
\newcommand{\EpiCommaCat}{\CommaCat[epi]}

\newcommand{\DaggerEpiCommaCat}{\CommaCat[dagger epi]}
\newcommand{\DaggerMonoCommaCat}{\CommaCat[dagger mono]}

\reversemarginpar

\begin{document}
\title[Dagger categories and the complex numbers]{Dagger categories and the complex numbers\\[1.5ex]\footnotesize Axioms for the category of finite-dimensional\\Hilbert spaces and linear contractions}

\date{\today}
\author{Matthew {Di Meglio}}
\author{Chris Heunen}
\address{University of Edinburgh}
\email{m.dimeglio@ed.ac.uk}
\email{chris.heunen@ed.ac.uk}

\begin{abstract}
We unravel a deep connection between limits of real numbers and limits in category theory. Using a new variant of the classical characterisation of the real numbers, we characterise the category of finite-dimensional Hilbert spaces and linear contractions in terms of simple category-theoretic structures and properties that do not refer to norms, continuity, or real numbers. This builds on Heunen, Kornell, and Van der Schaaf's easier characterisation of the category of all Hilbert spaces and linear contractions.
\end{abstract}
\maketitle

\section{Introduction}

The category \(\Hilb\) of Hilbert spaces and bounded linear maps and the category \(\Con\) of Hilbert spaces and linear contractions were both recently characterised in terms of simple category-theoretic structures and properties~\cite{heunenkornell:hilb,heunenkornellvanderschaaf:con}. For example, the structure of a \emph{dagger} encodes adjoints of linear maps. Remarkably, none of these properties refer to analytic notions such as norms, continuity, dimension, real numbers, convexity or probability. For mathematicians, these characterisations give a surprisingly new perspective on Hilbert spaces—a well-studied structure in functional analysis. For theoretical physicists, they provide further justification for the category-theoretic approach to quantum mechanics~\cite{heunenvicary:cqm}.

In quantum computing and quantum information theory, the Hilbert spaces of interest are typically finite dimensional.
Counterintuitively, finding axioms for categories with only \textit{finite-dimensional} Hilbert spaces is more challenging than doing so for categories with \textit{all} Hilbert spaces. The issue is that the natural category-theoretic way to encode analytic completeness of the scalar field is in terms of directed colimits, but the existence of too many of these colimits also implies the existence of objects corresponding to infinite-dimensional spaces. Until now, the only known way to prove that the scalars are the real or complex numbers was to construct such an infinite-dimensional object and then apply Solèr's theorem~\cite{soler}. Without such infinite-dimensional objects, a different approach is necessary.

An obvious way to bypass Solèr's theorem is to directly appeal to the classical characterisation of the reals as the unique Dedekind-complete ordered field, but it is unclear how to prove that the scalars have these specific properties. \citeauthor{demarr:partiallyorderedfields} showed that the reals are also the unique partially ordered field with suprema of bounded increasing sequences~\cite{demarr:partiallyorderedfields}. 
Defining and ordering \textit{positive} scalars 
based on the observations that each positive real is the squared norm of some vector and that contractions decrease norms, the supremum of a bounded increasing sequence of positive scalars can be recovered from the colimit of a directed diagram associated to the sequence. This explicit construction of limits in real analysis from limits in category theory, together with an extension of \citeauthor{demarr:partiallyorderedfields}'s theorem for partially ordered semifields that embed nicely in a field, is our first contribution.

Our second contribution is a resolution of the tension between too few and too many directed diagrams having colimits. We identify the \textit{bounded} ones—the ones that admit a cocone of monomorphisms—as striking this balance. Colimits of these diagrams suffice to construct suprema as explained above. Yet the class of finite-dimensional Hilbert spaces is closed under taking these colimits, essentially because the domain of a monomorphism has dimension at most that of its codomain.

The final ingredient is a category-theoretic property enforcing \textit{finite dimensionality}. 
For this, we adopt a notion from operator algebra~\cite{murrayvonneumann:ringsofoperators}, which we call \emph{dagger finiteness}, that is closely related to Dedekind finiteness from set theory~\cite{ferreiros:history}, and Hopfianness from group theory~\cite{varadarajan:hopfian}. It is defined purely in terms of the dagger and composition. 

Combining these ideas, we give an axiomatic characterisation of the category \(\FCon\) of finite-dimensional Hilbert spaces and linear contractions. The axioms are listed in \cref{sec:axioms}. Most are identical to ones~\cite{heunenkornellvanderschaaf:con} for \(\Con\), so we keep our discussion of them brief. The high-level structure of our proof is also the same: show that the scalar localisation of a category satisfying our axioms is equivalent to the category of finite-dimensional Hilbert spaces and \textit{all} linear maps, then identify the original category with the full subcategory of linear contractions. Much of the proof~\cite{heunenkornellvanderschaaf:con} of the characterisation of \(\Con\) depends only on axioms that we retain. In \cref{sec:prelims} we recall the results that we reuse. Our proof then proceeds as follows: 
\begin{itemize}
    \item in \cref{sec:order}, we construct the partially ordered semifield of positive scalars and show that it has suprema of bounded increasing sequences;
    \item in \cref{sec:limits}, we characterise the real and complex numbers among involutive fields with a partially ordered subsemifield of \textit{positive} elements, and use this characterisation to deduce that the scalars are the real or complex numbers;
    \item in \cref{sec:finite}, we show that the inner-product space associated to each object is finite dimensional;
    \item in \cref{sec:con}, we complete the characterisation of $\FCon$, and outline how to eliminate the use of Solèr's theorem from the characterisation~\cite{heunenkornellvanderschaaf:con} of $\Con$.
\end{itemize}
In \cref{sec:colimits:epi}, assuming a different completeness axiom, we prove that the scalars are again the real or complex numbers, this time using a new characterisation of the positive reals among partially ordered semifields. Establishing finite dimensionality from this alternative axiom remains an open problem.

Future work will characterise the category $\FHilb$ of finite-dimensional Hilbert spaces and \emph{all} linear maps using similar ideas, also accounting for quaternionic Hilbert spaces by removing the need for a tensor product altogether.

\section{Axioms}
\label{sec:axioms}
A \textit{dagger} is a contravariant endofunctor \(\blank^\dagger\) such that $X^\dag = X$ for all objects~$X$ and $f^{\dag \dag} = f$ for all morphisms~$f$. A \textit{dagger category} is a category with a dagger. A morphism \(f\) is called a \emph{dagger monomorphism}
if $f^\dag f = \id$, a \textit{dagger epimorphism} if \(ff^\dagger = 1\), and a \textit{dagger isomorphism} if it is both a dagger monomorphism and a dagger epimorphism. Similar to how monomorphisms and epimorphisms \(X \to Y\) are often drawn as \(\mono{X}{Y}\) and \(\epi{X}{Y}\), respectively, dagger monomorphisms and dagger epimorphisms \(X \to Y\) will be drawn as \(\daggerMono{X}{Y}\) and \(\daggerEpi{X}{Y}\), respectively. A functor \(F\) between dagger categories is a \emph{dagger functor} if \(F(f^\dagger) = F(f)^\dagger\) for all morphisms \(f\). An \emph{equivalence} of dagger categories is a dagger functor that is full, faithful, and dagger essentially surjective, that is, every object in its codomain is dagger isomorphic to an object in its image.

A \emph{dagger symmetric monoidal category} is a symmetric monoidal category with a dagger whose monoidal product is a dagger functor and whose symmetry, associator and unitors are dagger isomorphisms. A strong monoidal functor between dagger symmetric monoidal categories is \emph{dagger strong monoidal} if it is a dagger functor and its coherence natural transformations are dagger isomorphisms. An \emph{equivalence} of dagger symmetric monoidal categories is a dagger strong monoidal functor that is also an equivalence of dagger categories.

A \textit{dagger rig category} is a dagger category equipped with dagger symmetric monoidal structures \((\otimes, I)\) and \((\oplus, O)\), and natural dagger isomorphisms
\begin{align*}
A \otimes (X \oplus Y) & \longrightarrow (A \otimes X) \oplus (A \otimes Y), &
X \otimes O & \longrightarrow O, \\
(X \oplus Y) \otimes A & \longrightarrow (X \otimes A) \oplus (Y \otimes A), &
O \otimes X & \longrightarrow O, 
\end{align*}
subject to certain coherence conditions~\cite[Section~1]{laplaza:rig}. An \emph{equivalence} of dagger rig categories is an equivalence of the underlying dagger categories whose underlying dagger functor is equipped with dagger strong monoidal structures for $\otimes$ and $\oplus$.

The goal of this article is to prove the following theorem. 

\begin{theorem}\label{thm:main}
  A locally small dagger rig category is equivalent to the dagger rig category $\FCon$ of finite-dimensional Hilbert spaces and linear contractions if and only if it satisfies \cref{axiom:affine,axiom:jointlyepic,axiom:nondegenerate,axiom:simple,axiom:separator,axiom:equalisers,axiom:kernels,axiom:positive,axiom:colimits:mono,axiom:finite} below.
\end{theorem}

\subsection{Familiar axioms} \cref{axiom:affine,axiom:jointlyepic,axiom:nondegenerate,axiom:simple,axiom:separator,axiom:equalisers,axiom:kernels,axiom:positive} are all also axioms for \(\Con\)~\cite[Section~2]{heunenkornellvanderschaaf:con}, so we keep our discussion of them here brief. The proof that $\Con$ satisfies these axioms~\cite[Section~3]{heunenkornellvanderschaaf:con} works \textit{mutatis mutandis} for $\FCon$.

\begin{axiom}
\label{axiom:affine}
The monoidal structure \((\oplus, O)\) is \textit{semicartesian} (or \textit{affine}).
\end{axiom}

This means that the object~$O$ is terminal. The dagger then makes it a zero object. Denote zero morphisms by \(0\), the injections
\[
    \begin{tikzcd}[cramped]
        X \cong X \oplus O \arrow[r, "1 \oplus 0"] \&
        X \oplus Y
    \end{tikzcd}
    \qquad\text{and}\qquad
    \begin{tikzcd}[cramped]
        Y \cong O \oplus Y \arrow[r, "0 \oplus \id"] \&
        X \oplus Y
    \end{tikzcd}
\]
by \(\inl\) and \(\inr\), respectively, and let \(p_1 = {i_1}^\dagger\) and \(p_2 = {i_2}^\dagger\).

\begin{axiom}
\label{axiom:jointlyepic}
The injections \(\inl \colon X \to X \oplus Y\) and \(\inr \colon Y \to X \oplus Y\) are \emph{jointly epic}.  
\end{axiom}

This means if \(f\inl = g\inl\) and \(f\inr = g\inr\), then \(f = g\).

\begin{axiom}
\label{axiom:nondegenerate}There is a morphism $d \colon I \to I \oplus I$ such that $\inl^\dag d \neq 0 \neq \inr^\dag d$.
\end{axiom}

\begin{axiom}
\label{axiom:simple}
The object $I$ is \emph{dagger simple}.
\end{axiom}

This means there are exactly two subobjects of \(I\) that have a dagger monic representative. These are necessarily $0 \colon \mono{O}{I}$ and $1 \colon \mono{I}{I}$.

\begin{axiom}
\label{axiom:separator}
The object $I$ is a \emph{\(\otimes\)-monoidal separator}.
\end{axiom}

This means if $f (x \otimes y) = g (x \otimes y)$ for all $x \colon I \to X$ and $y \colon I \to Y$, then \(f = g\).

\begin{axiom}
\label{axiom:equalisers}
Every parallel pair of morphisms has a dagger equaliser.
\end{axiom}

A \textit{dagger equaliser} is an equaliser that is a dagger monomorphism.

\begin{axiom}
\label{axiom:kernels}
Every dagger monomorphism is a kernel.
\end{axiom}

\begin{axiom}
\label{axiom:positive}
For all epimorphisms \(x \colon \epi{A}{X}\) and \(y \colon \epi{A}{Y}\), we have \(x^\dag x = y^\dag y\) if and only if there is an isomorphism \(f \colon X \to Y\) such that \(y = fx\).
\end{axiom}

\subsection{Completeness axiom}
All characterisations of the real numbers involve an infinitary assumption, such as Dedekind or Cauchy completeness. To ensure that the field of scalars is \(\Reals\) or \(\Comps\), we thus also need an infinitary axiom. The one used in the characterisation~\cite[Section~2]{heunenkornellvanderschaaf:con} of \(\Con\)—that every directed diagram has a colimit—does not hold in \(\FCon\). For example, the directed diagram
\[
\begin{tikzcd}[cramped, sep=large]
\Comps \arrow[r, "i_1"] \&
\Comps^2 \arrow[r, "i_{1,2}"] \&
\Comps^3 \arrow[r, "i_{1,2,3}"]\&
\ldots,
\end{tikzcd}
\]
whose colimit in \(\Con\) is infinite dimensional, does not have a colimit in \(\FCon\). We will instead use a weakening of this assumption, which may also be viewed as a categorification of the condition~\cite{demarr:partiallyorderedfields} on a partially ordered field that every bounded increasing sequence has a supremum.

A \textit{sequential diagram} is a special kind of directed diagram; in particular, it is a diagram generated by a sequence of objects and morphisms of the form
\[
\begin{tikzcd}[cramped, sep=large]
    X_1
        \arrow[r, "f_1"]
        \&
    X_2
        \arrow[r, "f_2"]
        \&
    X_3
        \arrow[r, "f_3"]
        \&
    \cdots,
\end{tikzcd}
\]
which we sometimes abbreviate to \((X_n, f_n)\).\footnote{Sequential diagrams are precisely \(\omega\)-shaped diagrams, where \(\omega\) is the smallest limit ordinal.} Dually, a \textit{cosequential diagram} is a diagram generated by a sequence of objects and morphisms of the form
\[
\begin{tikzcd}[cramped, sep=large]
    \cdots
        \arrow[r, "f_3"]
        \&
    X_3
        \arrow[r, "f_2"]
        \&
    X_2
        \arrow[r, "f_1"]
        \&
    X_1,
\end{tikzcd}   
\]
which we sometimes also abbreviate to \((X_n, f_n)\).

Call a diagram \emph{bounded} when it has a cocone of monomorphisms and \textit{cobounded} when it has a cone of epimorphisms. The morphisms in a bounded diagram are necessarily monic whilst those in a cobounded diagram are necessarily epic.

\begin{axiom}
\label{axiom:colimits:mono}
Every bounded sequential diagram has a colimit.
\end{axiom}

With the dagger, this is equivalent to every cobounded cosequential diagram having a limit. We will swap between considering colimits of bounded sequential diagrams and limits of cobounded cosequential diagrams as is convenient.

\begin{lemma}
The category \(\FCon\) satisfies \cref{axiom:colimits:mono}.
\end{lemma}

\begin{proof}
Let $(X_n,f_n)$ be a bounded sequential diagram in $\FCon$, and let $a_n \colon X_n \to A$ be a cocone of monomorphisms. Let $C$ be the union in \(A\) of the images of the \(a_n\) and define $c_n \colon X_n \to C$ by $c_nx=a_nx$. As \(C\) is a vector subspace of \(A\), which is a finite-dimensional Hilbert space, the restriction of the inner product of \(A\) to \(C\) makes \(C\) into another finite-dimensional Hilbert space. The maps \(c_n\) then form a cocone on the diagram \((X_n,f_n)\) in \(\FCon\). To see that it is an initial cocone, let $b_n \colon X_n \to B$ be another cocone. Define \(m \colon C \to B\) by \(m a_n x = b_n x\). This map is well defined because \(b_n\) is a cocone. It is also total because every element of \(C\) is of the form \(a_n x\) for some \(n\) and some \(x \in X_n\). It is actually a linear contraction, and, in particular, the unique one satisfying \(m c_n = b_n\).
\end{proof}

\subsection{Finiteness axiom} Both \(\Con\) and \(\FCon\) satisfy all of the axioms listed so far. Distinguishing between these two categories requires an axiom that encodes finite dimensionality. The notion of \textit{dagger finiteness}, defined below, comes from operator algebra~\cite[Section~7.1]{murrayvonneumann:ringsofoperators}. It is similar to the notion of \textit{Dedekind finiteness} from set theory~\cite[Section~5.2.2]{ferreiros:history}, which has also been adapted to other types of categories (see~\cite[Theorem~1.1]{leary:dedekindfinite} and ~\cite[Definition~1.1]{stout:dedekindfinite}).

\begin{definition}\label{def:daggerfinite}
An object \(X\) is called \textit{dagger finite} when, for each \(f \colon X \to X\), if \(f^\dagger f = 1\) then \(ff^\dagger = 1\).
\end{definition}

In other words, an object \(X\) is dagger finite if every dagger monic endomorphism on \(X\) is a dagger isomorphism.

\begin{axiom}
\label{axiom:finite}
Every object is dagger finite.
\end{axiom}

We now show that this axiom holds in \(\FCon\) and not in \(\Con\).

\begin{proposition}
\label{prop:con_dagger_finite}
A Hilbert space is dagger finite in \(\Con\) if and only if it is finite dimensional. In particular, every object of \(\FCon\) is dagger finite.
\end{proposition}

\begin{proof}
Let \(X\) be a Hilbert space. If \(X\) is finite dimensional, then \(X\) is dagger finite by the rank-nullity theorem. Conversely, suppose that \(X\) is infinite dimensional. Then it contains as a closed subspace a copy of the Hilbert space \(\ell_2(\Nats)\) of square summable sequences. The right shift map, which sends the $n$th standard basis vector to the $(n+1)$th standard basis vector, is a dagger monic endomorphism on \(\ell_2(\Nats)\) that is not a dagger isomorphism. Pairing it with the identity map on the orthogonal complement of \(\ell_2(\Nats)\) in \(X\), we obtain a dagger monic endomorphism on~\(X\) that is not a dagger isomorphism. Hence \(X\) is not dagger finite. 
\end{proof}

Similarly, an object of \(\Hilb\) is dagger finite if and only if it is finite dimensional. Also, every object of \(\FHilb\) is dagger finite. In \cref{sec:finite}, we will use an abstract version of this argument to prove that every dagger finite object \(X\) in a dagger rig category satisfying \cref{axiom:affine,axiom:jointlyepic,axiom:nondegenerate,axiom:simple,axiom:separator,axiom:equalisers,axiom:kernels,axiom:positive,axiom:colimits:mono} is dagger isomorphic to \(I \oplus I \oplus \dots \oplus I\). It will follow that the inner-product space corresponding to \(X\) is finite dimensional.

\section{The scalar localisation}
\label{sec:prelims}

Let \(\D\) be a locally small dagger rig category that satisfies \cref{axiom:affine,axiom:jointlyepic,axiom:nondegenerate,axiom:simple,axiom:separator,axiom:equalisers,axiom:kernels,axiom:positive,axiom:colimits:mono,axiom:finite}. Our goal for the remainder of the article is to prove that \(\D\) is equivalent to \(\FCon\). We begin by recalling some constructions and results that follow from \cref{axiom:affine,axiom:jointlyepic,axiom:nondegenerate,axiom:simple,axiom:separator,axiom:equalisers,axiom:kernels,axiom:positive}.

The set 
\[\Disk = \set[\big]{a \colon I \to I \text{ in } \D}\]
of \emph{scalars} of \(\D\) is a commutative absorption monoid under composition, with unit the identity morphism \(1 \colon I \to I\) and absorbing element the zero morphism \(0 \colon I \to I\). If \(\D\) is \(\Con\) or \(\FCon\), then \(\Disk\) is the unit disk in \(\Reals\) or \(\Comps\).

The monoidal structure \(\otimes\) on \(\D\) induces an action of the absorption monoid \(\Disk\) on the category \(\D\) called \emph{scalar multiplication}. This action is defined by the equation
\[
    a \cdot f
    = \paren[\big]{\begin{tikzcd}[cramped]
        X
            \arrow[r, "\lambda^{-1}"]
            \&
        I \otimes X
            \arrow[r, "a \otimes f"]
            \&
        I \otimes Y
            \arrow[r, "\lambda"]
            \&
        Y
    \end{tikzcd}}
\]
for each scalar \(a \in \Disk\) and each morphism \(f \colon X \to Y\) of~\(\D\).

The \textit{dagger monoidal localisation} of a dagger monoidal category at a class of morphisms, if it exists, is the initial dagger strong monoidal functor out of the category that sends all morphisms in the class to isomorphisms. We are interested in the dagger \(\otimes\)-monoidal localisation \(U \colon \D \to \C\) of \(\D\) at the set
\[\NZDisk = \setb{a \in \Disk}{a \neq 0}.\]
It has the following concrete description~\cite[Proposition~10, Lemmas~12–14]{heunenkornellvanderschaaf:con}. The objects of \(\C\) are the same as the objects of \(\D\). The morphisms \(X \to Y\) of \(\C\) are the equivalence classes of the equivalence relation \(\simeq\) on
\[\setb[\big]{(f, a)}{f \colon X \to Y \text{ in } \D, \; a \in \NZDisk}\]
defined by 
\[
    (f, a) \simeq (g, b)
    \quad\iff\quad
    b \cdot f = a \cdot g.
\]
The equivalence class of \((f, a)\) will be represented by the fraction \(f/a\). The dagger on \(\C\) is defined by \((f/a)^\dagger = f^\dagger/a^\dagger\). The action of \(\otimes\) on the objects of \(\C\) is the same as its action on the objects of \(\D\), and its action on the morphisms of \(\C\) is defined by
\[
    \frac{f}{a} \otimes \frac{g}{b} = \frac{f \otimes g}{ab}.
\]
The functor \(U \colon \D \to \C\) is the identity on objects and is defined by \(f \mapsto f/1\) on morphisms. It exhibits \(\D\) as a dagger \(\otimes\)-monoidal wide subcategory of \(\C\). We will often identify morphisms of \(\D\) with their \(U\)-image in \(\C\), writing \(f\) instead of \(f/1\).

There is also a \(\oplus\)-monoidal structure on \(\C\) that makes \(U\) strict \(\oplus\)-monoidal. The action of \(\oplus\) on the objects of \(\C\) is the same as its action on the objects of \(\D\), and its action on the morphisms of \(\C\) is defined by
\[
  \frac{f}{a} \oplus \frac{g}{b} = \frac{b \cdot f \oplus a \cdot g}{ab}.
\]
Actually~\cite[Lemma~17]{heunenkornellvanderschaaf:con}, \(O\) is a zero object in \(\C\), and, for all objects \(X\) and \(Y\),
\[
    \begin{tikzcd}
        X 
            \arrow[r, "i_1", shift left]
            \arrow[from=r, "p_1", shift left] \&
        X \oplus Y \&
        Y
            \arrow[l, "i_2" swap, shift right]
            \arrow[from=l, "p_2" swap, shift right]
    \end{tikzcd}
\]
is a dagger biproduct in \(\C\). We use the usual matrix notation for morphisms between these biproducts. For example, given morphisms \(f_{jk} \colon X_k \to Y_j\) for all \(j,k \in \set{1, 2}\), the matrix \(\begin{bsmallmatrix} f_{11} & f_{12} \\ f_{21} & f_{22} \end{bsmallmatrix}\) represents the unique morphism \(f \colon X_1 \oplus X_2 \to Y_1 \oplus Y_2\) such that \(p_j f i_k = f_{jk}\) for all \(j, k \in \set{1,2}\).

The category \(\C\) has the following additional properties: every parallel pair of morphisms has a dagger equaliser~\cite[Lemma~18]{heunenkornellvanderschaaf:con}, every dagger monomorphism is a kernel~\cite[Lemma~19]{heunenkornellvanderschaaf:con}, the object \(I\) is simple~\cite[Lemma~16]{heunenkornellvanderschaaf:con}, and it is also a \(\otimes\)-monoidal separator~\cite[Lemma~17]{heunenkornellvanderschaaf:con}.
The semiring 
\[
  \Scalars = \set[\big]{a \colon I \to I \text{ in } \C}
\]
of scalars in \(\C\) is thus an involutive field of characteristic zero~\cite[Theorem~4.8]{heunen:hilbert-categories}. Let \(\SAScalars\) be its subfield of self-adjoint scalars. If \(\D\) is \(\FCon\), then \(\C\) is \(\FHilb\), and \(\SAScalars\) is \(\Reals\), and \(\Scalars\) is \(\Reals\) or \(\Comps\).

\section{The partially ordered semifield of positive scalars}\label{sec:order}

In this section, we distinguish a set \(\PosScalars\) of \emph{positive}\footnote{In this article, the terms \textit{positive} and \textit{negative} include zero, and the terms \textit{increasing} and \textit{decreasing} include equality.
} scalars of \(\C\) and equip it with the structure of a \textit{partially ordered strict semifield}. We recall the precise definition in due course. For now, it suffices to know that a strict semifield is like a field but without additive inverses, and a partially ordered strict semifield is a strict semifield equipped with a partial order that appropriately respects the semifield operations.

The positive scalars should correspond to the squared norms of vectors, so that one positive scalar is larger than another exactly when there is a contraction that maps a vector representing the first to a vector representing the second. Abstractly, these vectors and contractions are the objects and morphisms of the comma category~\(\CommaCat{I}{U}\). Concretely, its objects are pairs \((X, x)\) where \(x \colon I \to X\) is a morphism in \(\C\), and its morphisms \(f \colon (X, x) \to (Y, y)\) are the morphisms \(f \colon X \to Y\) in \(\D\) such that \(y = fx\). Let \(\normSq\) be the function from the set of objects of \(\CommaCat{I}{U}\) to \(\Scalars\) that maps \((X, x)\) to \(x^\dag x\), and let \(\PosScalars\) be the image of \(\normSq\). The elements of \(\PosScalars\) will be called \textit{positive} scalars.

Our initial goal is to define a partial order on~\(\PosScalars\) so that \(\normSq \colon \CommaCat{I}{U} \to \PosScalars\) is functorial.
Every category \(\A\) has a universal collapse to a partially ordered class,
namely its \textit{partially ordered reflection} \(\Par \A\), which is described concretely below. By \cref{Prop: Preorder isomorphism and norms} below, the function \(\normSq\) factors through the object map of the canonical functor \(\CommaCat{I}{U} \to \Par (\CommaCat{I}{U})\) via a bijection, and this bijection induces the desired partial order on~\(\PosScalars\).

Concretely, the elements of \(\Par \A\) are the equivalence classes of the equivalence relation \(\simeq\) on the class of objects of \(\A\) defined by \(A \simeq B\) if there exist morphisms \(A \to B\) and \(B \to A\).  Write \([A]\) for the equivalence class of an object \(A\) of \(\A\). The partial order \(\geq\) of \(\Par \A\) is then defined by \([A] \geq [B]\) if there is a morphism \(A \to B\). The canonical functor \(\A \to \Par \A\) maps each object \(A\) to its equivalence class \([A]\), and is uniquely determined on morphisms.

\begin{proposition}
\label{Prop: Preorder isomorphism and norms}
Let \((X,x)\) and \((Y,y)\) be objects of \(\CommaCat{I}{U}\). Then \((X,x) \simeq (Y,y)\) if and only if \(x^\dag x = y^\dag y\).
\end{proposition}

\Cref{Prop: Preorder isomorphism and norms} is really the analogue of \cref{axiom:positive} for \(\C\). To prove it, we need the following two lemmas. The first allows us to focus on those objects \((X, x)\) of \(\CommaCat{I}{U}\) where \(x\) is epic in \(\C\). The second relates the epimorphisms in~\(\C\) and those in \(\D\).

\begin{lemma}
\label{lem:epicrepresentative}
  For each object \((X, x)\) of \(\CommaCat{I}{U}\), there is an object \((E, e)\) of \(\CommaCat{I}{U}\) with \(e\) epic in \(\C\) such that \((E, e) \simeq (X,x)\) and \(e^\dag e = x^\dag x\).
\end{lemma}

\begin{proof}
  As \(\C\) has (epic, dagger monic) factorisations, there is an epimorphism \(e \colon I \to E\) and a dagger monomorphism \(m \colon E \to X\) such that \(x = me\). The morphism \(m\) comes from \(\D\)~\cite[Lemma~14]{heunenkornellvanderschaaf:con}. As \(\D\) is a \textit{dagger} subcategory of \(\C\), so does \(m^\dagger\). Hence \([(E, e)] \geq [(X, x)]\). Also \([(X, x)] \geq [(E, e)]\) because \(m^\dagger x = m^\dagger m e= e\). Finally, we have \(x^\dagger x = e^\dagger m^\dagger m e = e^\dagger e\).
\end{proof}

\begin{lemma}
\label{lem:monos}
The embedding \(U \colon \D \to \C\) preserves and reflects epimorphisms.
\end{lemma}

\begin{proof}
  Reflection follows from the faithfulness of \(U\). For preservation, let \(e \colon \epi{A}{X}\) be an epimorphism in \(\D\). Let \(f, g \colon X \to Y\) in \(\C\), and suppose that \(fe = ge\). Now \(f = s/a\) and \(g = t/b\) for some \(s, t \colon X \to Y\) in \(\D\) and some \(a, b \in \NZDisk\). As
  \[
    \frac{se}{a} = \frac{s}{a} e = fe = ge = \frac{t}{b} e = \frac{te}{b},
  \]
  we have \((b \cdot s)e = b \cdot se = a \cdot te = (a \cdot t) e\). As \(e\) is epic in \(\D\), actually \(b \cdot s = a \cdot t\), and so \(f = s/a = t/b = g\).
\end{proof}

\begin{proof}[{Proof of \cref{Prop: Preorder isomorphism and norms}}]
By \cref{lem:epicrepresentative}, we may assume, without loss of generality, that \(x\) and \(y\) are epic in \(\C\). Now \(x = u/a\) and \(y = v/b\) for some morphisms \(u \colon I \to X\) and \(v \colon I \to Y\) in \(\D\) and some scalars \(a, b \in \NZDisk\). As \(x\) and \(y\) are epic in \(\C\) and \(a\) and \(b\) are invertible in \(\C\), their composites \(u = xa\) and \(v = yb\) are epic in \(\C\), and thus, by \cref{lem:monos}, also epic in \(\D\). Cross-multiplying, we see that \(x^\dagger x  = y^\dagger y\) if and only if \((ub)^\dagger ub = (va)^\dagger va\). By \cref{axiom:positive}, the latter equation holds exactly when there is an isomorphism \(f \colon X \to Y\) in \(\D\) such that \(va = fub\) in \(\D\), or equivalently, such that \(y = fx\) in \(\C\). If such an isomorphism \(f\) exists, clearly \((X,x) \simeq (Y,y)\). Conversely, suppose that \((X, x) \simeq (Y, y)\). Then there are morphisms \(f \colon X \to Y\) and \(g \colon Y \to X\) in \(\D\) such that \(y = fx\) and \(x = gy\). As \(fgy = fx = y\) and \(y\) is epic, actually \(fg = 1\). Similarly \(gf = 1\), so \(f\) is actually an isomorphism.
\end{proof}

\begin{lemma}\label{lem:monosinD}
  A morphism \(x \colon I \to X\) in \(\C\) comes from \(\D\) exactly when \(1 \geq x^\dag x\).
\end{lemma}
\begin{proof}
  This is an exercise in unpacking definitions. The inequality \(1 \geq x^\dag x\) holds if and only if there exists a morphism \(f \colon (I, 1) \to (X, x)\) in \(\CommaCat{I}{U}\), that is, if and only if there exists a morphism \(f \colon I \to X\) in \(\D\) such that \(x = f1\).
\end{proof}

A \textit{semifield} is a set \(S\) equipped with two binary operations~\(+\) and~\(\cdot\), called \textit{addition} and \textit{multiplication}, and two distinct distinguished elements~\(0\) and~\(1\), called \textit{zero} and \textit{one}, such that \((S, +, 0)\) and \((S, \cdot, 1)\) are commutative monoids, multiplication distributes over addition, every non-zero element has a multiplicative inverse, and every element is annihilated by zero. Every field is a semifield. A semifield is called \textit{strict} if it is not a field, or, equivalently, if \(1\) does not have an additive inverse.

\begin{proposition}
The set \(\PosScalars\) is a subsemifield of \(\Scalars\).
\end{proposition}

\begin{proof}
We have \(0 = 0^\dag 0\) and \(1 = 1^\dag 1\), so \(\PosScalars\) contains \(0\) and \(1\). For all morphisms \(x \colon I \to X\) and \(y \colon I \to Y\) in \(\C\), we have
\begin{align*}
    x^\dag x + y^\dag y &= \Delta^\dag (x \oplus y)^\dag (x \oplus y) \Delta,\\
    x^\dag x \cdot y^\dag y &= \lambda (x \otimes y)^\dag (x \otimes y) \lambda^\dag,\\
    \intertext{and, whenever \(x^\dag x \neq 0\), also}
    (x^\dag x)^{-1} &= (x^\dag x)^{-1\dag}x^\dag x(x^\dag x)^{-1},
\end{align*}
so \(\PosScalars\) is also closed under addition, multiplication, and inversion.
\end{proof}

\begin{proposition}
The semifield \(\PosScalars\) is strict.
\end{proposition}

\begin{proof}
Suppose that \(-1 \in \PosScalars\). Then \(-1 = x^\dag x\) for some \(x \colon I \to X\) in \(\C\). Then \(\pair{1}{x} = 0\) because
\(\pair{1}{x}^\dag \pair{1}{x} = 1 + x^\dag x = 0\)
and \(\C\) has dagger equalisers~\cite[Lemma~II.5]{vicary:complex}. Hence
\(1 = \pi_1 \pair{1}{x} = \pi_1 0 = 0,\)
which is a contradiction because \(\Scalars\) is a field.
\end{proof}

We will use the following variant of \citeauthor{fritz:semifield}'s \textit{preordered semifield}~\cite[Definition~3.16]{fritz:semifield}. Unlike \citeauthor{fritz:semifield}, we incorporate the assumption \(1 \geq 0\) into the definition instead of restating it every time we need it. 

\begin{definition}
A \textit{partially ordered strict semifield} is a strict semifield equipped with a partial order \(\geq\) satisfying the following axioms:
\begin{itemize}
    \item
    \textit{Addition is monotonic:} if \(a \geq b\) then \(a + c \geq b + c\).
    \item 
    \textit{Multiplication is monotonic:} if \(a \geq b\) then \(ac \geq bc\).
    \item 
    \textit{One is positive:} \(1 \geq 0\).
\end{itemize}
An \textit{ordered strict semifield} is a partially ordered strict semifield whose order is total.
\end{definition}

The axioms above are tailored for semifields that are strict. They imply that \(a = 1a \geq 0a = 0\) for all elements \(a\), whereas \(-1 \leq 0\) in all partially ordered fields.

Examples of ordered strict semifields include the \textit{rational semifield} \(\PosRats\), the \textit{real semifield} \(\PosReals\), and the \textit{tropical semifield} \(\TropReals\). The rational and real semifields are, respectively, the sets of positive rational and real numbers, with their usual addition, multiplication and ordering. The tropical semifield is also the set of positive real numbers with its usual multiplication and ordering, but with `addition' given by maximum rather than sum.

\begin{proposition}
\label{prop:orderedsemifield}
The semifield \(\PosScalars\) is a partially ordered strict semifield when it is equipped with the partial order that it inherits from \( \Par(\CommaCat{I}{U})\).
\end{proposition}

\begin{proof}
Let \(x \colon I \to X\), \(y \colon I \to Y\) and \(z \colon I \to Z\) in \(\C\), and suppose that \(x^\dag x \geq y^\dag y\), that is, that there is a morphism \(f \colon X \to Y\) in \(\D\) such that \(y = fx\). As
\[(y \oplus z) \Delta = (fx \oplus z) \Delta = (f \oplus 1)(x \oplus z) \Delta\]
and \(f \oplus 1\) is in \(\D\), we have \(x^\dag x + z^\dag z \geq y^\dag y + z^\dag z\), so addition is monotonic. As
\[
(y \otimes z) \lambda^\dag = (fx \otimes z) \lambda^\dag = (f \otimes 1)(x \otimes z) \lambda^\dag
\]
and \(f \otimes 1\) is in \(\D\), we have \(x^\dag x \cdot z^\dag z \geq y^\dag y \cdot z^\dag z\), so multiplication is also monotonic. Finally \(1 \geq 0\) because \(0 = 0 \circ 1\) and \(0\) is in \(\D\). 
\end{proof}

\begin{remark}
\label{rem:decategorification}
The functor \(\normSq \colon \CommaCat{I}{U} \to \PosScalars\) exhibits the partially ordered semifield~\(\PosScalars\) of positive scalars as a decategorification of the rig category \(\CommaCat{I}{U}\). Direct sums become addition, tensor products become multiplication, and morphisms become the partial order. This decategorification mirrors the analogy between the operators on a Hilbert space and the complex numbers.
\end{remark}

Our goal for the remainder of this section is to prove the following proposition.

\begin{proposition}
\label{prop:posscalarssups}
The partially ordered strict semifield \(\PosScalars\) has suprema of bounded increasing sequences and these are preserved by every endomap of the form \(a + \blank\).
\end{proposition}

The analogous statement about \(\D\) is that it has limits of cobounded cosequential diagrams and these are preserved by every endofunctor of the form \(X \oplus \blank\). The existence of these limits is the dual of \cref{axiom:colimits:mono}. Their preservation is the dual of the following proposition.

\begin{proposition}
\label{lem:biproductspreservecolimits}
Every endofunctor on \(\D\) of the form \(X \oplus \blank\) preserves colimits of bounded sequential diagrams.
\end{proposition}

To prove this proposition, we will use the following three lemmas.

\begin{lemma}
\label{lem:biproducts}
In a semiadditive category in which every split monomorphism is a normal monomorphism, the diagram
\[
\begin{tikzcd}[cramped, sep=huge]
A_1
    \arrow[r, "s_1", shift left]
    \arrow[from=r, "r_1", shift left]
    \&
A
    \arrow[from=r, "s_2" swap, shift right]
    \arrow[r, "r_2" swap, shift right]
    \&
A_2
\end{tikzcd}
\]
is a biproduct if and only if 
\[
r_1 s_1 = 1,
\qquad
r_2 s_2 = 1,
\qquad
r_2 = \coker(s_1),
\qquad
r_1 s_2 = 0.
\]
\end{lemma}
\begin{proof}
The \textit{only if} direction is well known. For the \textit{if} direction, the morphism \(\copair{s_1}{s_2} \colon A_1 \oplus A_2 \to A\)  is a section of \(\pair{r_1}{r_2} \colon A \to A_1 \oplus A_2\) because
\[
\begin{bmatrix}r_1 \\ r_2 \end{bmatrix} \begin{bmatrix} s_1 & s_2 \end{bmatrix} = \begin{bmatrix} r_1s_1 & r_1s_2 \\ r_2s_1 & r_2s_2 \end{bmatrix} = \begin{bmatrix} 1 & 0 \\ 0 & 1 \end{bmatrix} = 1.
\]
As it is then a normal monomorphism, it is an isomorphism if it has cokernel zero.

Suppose that \(f\copair{s_1}{s_2}=0\). Then $fs_1=0$ and so $f=f_2r_2$ for some morphism $f_2$. But then $0=fs_2=f_2r_2s_2=f_2$, and so $f=0r_2=0$. Hence $\coker\copair{s_1}{s_2}=0$.
\end{proof}

\begin{lemma}
\label{lem:dagger_mono_epi_dagger}
In a dagger category with finite dagger biproducts and dagger equalisers, if \(m \colon A \to X\) and \(e \colon X \to A\) satisfy \(em = 1\), then \(e = m^\dagger\) if and only if \(m\) is dagger monic and \(e\) is dagger epic.
\end{lemma}

\begin{proof}
The \textit{only if} direction is trivial. The \textit{if} direction follows from the equation \(m^\dagger m + ee^\dagger = 1 + 1 = m^\dagger e^\dagger + em\) because the dagger category has finite dagger biproducts and dagger equalisers~\cite[Lemma~II.9]{vicary:complex}.
\end{proof}

\begin{lemma}
\label{lem:normal_mono_D}
Every normal monomorphism in \(\D\) is dagger monic. Dually, every normal epimorphism in \(\D\) is dagger epic.
\end{lemma}

\begin{proof}
Let \(m \colon A \to X\) be a normal monomorphism in \(\D\). By \cref{axiom:kernels}, there is a morphism \(f \colon X \to Y\) in \(\D\) such that \(m\) is a kernel of \(f\). Let \(k \colon K \to X\) be a dagger kernel of \(f\). Then there is an isomorphism \(u \colon A \to K\) in \(\D\) such \(m = ku\). But \(u\) is a dagger isomorphism~\cite[Lemma~8]{heunenkornellvanderschaaf:con}. Hence \(m\), being the composite of two dagger monomorphisms, is itself dagger monic.
\end{proof}

\begin{proof}[Proof of \cref{lem:biproductspreservecolimits}]
Let $(Y_n,g_n)$ be a bounded sequential diagram in \(\D\), and let $c_n \colon Y_n \to \colim Y_n$ be a colimit of this diagram. As $(X \oplus Y_n, 1 \oplus g_n)$ is also a bounded sequential diagram, it has a colimit $d_n \colon X \oplus Y_n \to \colim (X \oplus Y_n)$.
\[
\begin{tikzcd}[sep=large]
X
    \arrow[r, "1" swap]
    \arrow[d, "i_1", shift left]
    \arrow[from=d, "p_1", shift left]
    \arrow[rrr, "1", bend left=15, shift left=2]
    \&
X
    \arrow[d, "i_1", shift left]
    \arrow[from=d, "p_1", shift left]
    \arrow[r, "1" swap]
    \arrow[rr, "1", bend left=7.5, shift left=1, shorten >=1.5ex]
    \&
\cdots
    \&[-2em]
X
    \arrow[d, "s_1", shift left]
    \arrow[from=d, "r_1", shift left]
    \\
X \oplus Y_1
    \arrow[r, "1 \oplus g_1"]
    \&
X \oplus Y_2
    \arrow[r, "1 \oplus g_2"]
    \arrow[from=d, "i_2" swap, shift right, near start]
    \arrow[d, "p_2" swap, shift right, near end]
    \&
\cdots
    \&
\colim (X \oplus Y_n)
    \arrow[from=lll, crossing over, bend right=15, "d_1" swap]
    \arrow[from=ll, bend right=7.5, "d_2"{swap}]
    \\
Y_1
    \arrow[r, "g_1"]
    \arrow[u, "i_2" swap, shift right]
    \arrow[from=u, "p_2" swap, shift right]
    \&
Y_2
    \arrow[r, "g_2"]
    \&
\cdots
    \&
\colim Y_n
    \arrow[from=lll, crossing over, bend right=15, "c_1" swap, shift right]
    \arrow[from=ll, bend right=7.5, "c_2"{swap}, shorten >=1ex, shift right]
    \arrow[u, "s_2" swap, shift right]
    \arrow[from=u, "r_2" swap, shift right]
\end{tikzcd}
\]
By universality of $d$, there is a unique morphism $r_2 \colon \colim (X \oplus Y_n) \to \colim Y_n$ such that $r_2 d_n = c_n p_2$ for each $n$. Also, letting $s_1 = d_1 i_1$, we have $s_1 = d_n i_1$ for each $n$. As colimits commute with colimits, the morphism $r_2$ is actually a cokernel of $s_1$ in \(\D\), and thus also in \(\C\)~\cite[Lemma~18]{heunenkornellvanderschaaf:con}. Similarly, there are unique morphisms $s_2 \colon \colim Y_n \to \colim (X \oplus Y_n)$ and $r_1 \colon \colim (X \oplus Y_n) \to X$ such that $p_1 = r_1 d_n$ and $d_n i_2 = s_2 c_n$ for each $n$, and $r_1$ is a cokernel of $s_2$ in \(\C\). Now $r_2 s_2 = 1$ because $r_2 s_2 c_n = r_2 d_n i_2 = c_n p_2 i_2 = c_n$ and the morphisms $c_n$ are jointly epic. Also $r_1 s_1 = r_1 d_1 i_1 = p_1 i_1 = 1$. By \cref{lem:biproducts}, the tuple \(\paren[\big]{\colim(X \oplus Y_n), s_1, s_2, r_1, r_2}\) is a biproduct of \(X\) and \(\colim Y_n\) in \(\C\). It follows that \(s_1\) and \(s_2\) are, respectively,  kernels in \(\C\) of \(r_1\) and \(r_2\). 
By \cref{lem:normal_mono_D}, the morphisms \(s_1\) and \(s_2\) are dagger monic, and the morphisms \(r_1\) and \(r_2\) are dagger epic. It follows, by \cref{lem:dagger_mono_epi_dagger}, that \(r_1 = {s_1}^\dagger\) and \(r_2 = {s_2}^\dagger\). Hence \(\begin{bsmallmatrix} r_1 \\ r_2 \end{bsmallmatrix} \colon \colim (X \oplus Y_n) \to X \oplus \colim Y_n\) is a dagger isomorphism, and so comes from \(\D\)~\cite[Lemma~14]{heunenkornellvanderschaaf:con}. Thus \(\begin{bsmallmatrix} r_1 \\ r_2 \end{bsmallmatrix} d_n \colon X \oplus Y_n \to X \oplus \colim Y_n\) is another colimit cocone on the diagram \((X \oplus Y_n, 1 \oplus g_n)\). Finally \(\begin{bsmallmatrix} r_1 \\ r_2 \end{bsmallmatrix} d_n = 1 \oplus c_n\) because \(p_1 \begin{bsmallmatrix} r_1 \\ r_2 \end{bsmallmatrix} d_n = r_1 d_n = p_1\) and \(p_2 \begin{bsmallmatrix} r_1 \\ r_2 \end{bsmallmatrix} d_n = r_2 d_n = c_n p_2\).
\end{proof}

To prove \cref{prop:posscalarssups} from these properties of \(\D\), we consider first the forgetful functor \(\Pi \colon \CommaCat{I}{U} \to \D\) and then the functor \(\normSq \colon \CommaCat{I}{U} \to \PosScalars\) defined in \cref{sec:order}.

\begin{proposition}
\label{prop:picreateslimits}
The functor \(\Pi\) creates limits of diagrams that have a cone.
\end{proposition}

\begin{proof}
Let \(f_r \colon (X_j, x_j) \to (X_k, x_k)\) be the value at \(r \colon j \to k\) of a diagram in \(\CommaCat{I}{U}\), and let \(t_j \colon (Y,y) \to (X_j, x_j)\) be a cone on this diagram. Suppose that the diagram \(f_{r} \colon X_j \to X_k\) in \(\D\) has a limit cone \(s_j \colon X \to X_j\) in \(\D\). Then there is a unique morphism \(t \colon Y \to X\) in \(\D\) such that \(t_j = s_jt\) for all \(j\). Also, similarly to \cref{lem:monos}, the morphisms \(s_j\) are jointly monic in \(\C\) because they are jointly monic in \(\D\).

First, we show that there exists a unique morphism \(x \colon I \to X\) in \(\C\) such that \(s_j \colon (X, x) \to (X_j, x_j)\) is a cone on \(f_r \colon (X_j, x_j) \to (X_k, x_k)\) in \(\CommaCat{I}{U}\). For uniqueness, observe that if \(x\) exists, then \(s_j x = x_j = t_jy = s_jty\) for all~\(j\), and so \(x = ty\) because the morphisms \(s_j\) are jointly monic in \(\C\). For existence, let \(x = ty\). For each \(j\), the morphism \(s_j\) is a morphism from \((X,x)\) to \((X_j, x_j)\) in~\(\CommaCat{I}{U}\) because \(s_j x = s_jty = t_j y = x_j\). Also \(s_j \colon (X,x) \to (X_j, x_j)\) is a cone on \(f_r \colon (X_j, x_j) \to (X_k, x_k)\) because \(s_j \colon X \to X_j\) is a cone on \(f_r \colon X_j \to X_k\).

We now show that \(s_j \colon (X, x) \to (X_j, x_j)\) is a limit of \(f_{r} \colon (X_j, x_j) \to (X_k, x_k)\). As the cone \(t_j \colon (Y,y) \to (X_j, x_j)\) is arbitrary, it suffices to show that there is a unique morphism \(t' \colon (Y,y) \to (X, x)\) in \(\CommaCat{I}{U}\) such that \(s_j t' = t_j\) for all \(j\). For uniqueness, if \(t'\) exists, then \(s_jt' = t_j = s_j t\) for all \(j\), and so \(t' = t\) because the morphisms \(s_j\) are jointly monic in \(\D\). For existence, as \(t_j = s_jt\) for all \(j\), it suffices to show that \(t\) is a morphism from \((Y,y)\) to \((X,x)\) in~\(\CommaCat{I}{U}\). But \(s_jty = t_j y = x_j = s_jx\) for all \(j\), and the morphisms \(s_j\) are jointly monic in \(\C\), so \(ty = x\).
\end{proof}

\begin{remark}
\label{rem:picreatesepis}
As \(\Pi\) creates pushouts~\cite[Theorem~3]{burstall:computational} (see also \cite[Proposition~3.3.8]{riehl:categorycontext}), it also creates epimorphisms. This means that a morphism \(e \colon (X,x) \to (Y,y)\) in \(\CommaCat{I}{U}\) is epic if and only if the morphism \(e \colon X \to Y\) in \(\D\) is epic. 
\end{remark}

For each object \((X,x)\) of \(\CommaCat{I}{U}\), let \((X,x) \oplus \blank \colon \CommaCat{I}{U} \to \CommaCat{I}{U}\) be the functor that maps \(f \colon (A,a) \to (B,b)\) to \(1 \oplus f \colon (X \oplus A, \pair{x}{a}) \to (X \oplus B, \pair{x}{b})\).

\begin{corollary}
\label{prop:commacatlimits}
The category \(\CommaCat{I}{U}\) has limits of cobounded cosequential diagrams and these are preserved by every endofunctor of the form \((X,x) \oplus \blank\).
\end{corollary}

\begin{proof}
Existence follows from \cref{axiom:colimits:mono,prop:picreateslimits,rem:picreatesepis}. As
\[\begin{tikzcd}[column sep=large]
\CommaCat{I}{U}
  \arrow[d, "{(X,x) \oplus
  \blank}" swap]
  \arrow[r, "\Pi"]
  \&
\D
  \arrow[d, "{X \oplus \blank}"]
  \\
\CommaCat{I}{U}
  \arrow[r, "\Pi" swap]
  \&
\D
\end{tikzcd}\]
commutes, preservation follows from \cref{prop:picreateslimits,lem:biproductspreservecolimits,rem:picreatesepis}.
\end{proof}

Let $\EpiCommaCat{I}{U}$ be the full subcategory of \(\CommaCat{I}{U}\) spanned by the objects \((X,x)\) where \(x\) is epic in \(\C\), and let \(J\) denote the canonical embedding \(\EpiCommaCat{I}{U} \hookrightarrow \CommaCat{I}{U}\). For each morphism \(f \colon (X,x) \to (Y,y)\) in \(\EpiCommaCat{I}{U}\), the morphism \(f \colon X \to Y\) is epic in \(\C\) by epimorphism cancellation, and thus epic in \(\D\) by \cref{lem:monos}.

\begin{remark}
\label{lem:normSq:factors}
The functor \(\normSq \circ J\) is surjective on objects and full by \cref{lem:epicrepresentative}, and its domain category is thin.\footnote{It follows that \(\normSq \circ J\) is an equivalence, but this does not appear to simplify what follows.}
Hence all diagrams in \(\PosScalars\) factor through \(\normSq \circ J\), as do all cones and cocones on such diagrams.
\end{remark}

\begin{lemma}
\label{lem:factornormsq}
Every bounded increasing sequence in \(\PosScalars\) factors through \(\normSq\) via a cobounded cosequential diagram in \(\CommaCat{I}{U}\).  
\end{lemma}

\begin{proof}
Regarding \(\PosScalars\) as a category, increasing sequences and their upper bounds are precisely cosequential diagrams and their cones. These factor through \(\normSq \circ J\) as described in \cref{lem:normSq:factors}. Composing with \(J\) then gives a factorisation of them through \(\normSq\). The components of the resulting cones are epic by \cref{rem:picreatesepis}.
\end{proof}

\begin{lemma}
\label{lem:normsq_preserve}
The functor \(\normSq\) preserves limits.
\end{lemma}

\begin{proof}
Given a diagram \(D\) in \(\CommaCat{I}{U}\), every cone on \(\normSq \circ D\) is the image by \(\normSq\) of a cone on \(D\). To construct the lifted cone, \cref{lem:epicrepresentative} gives its apex, fullness of \(\normSq\) gives its components, and epicness of its apex gives its naturality.
\end{proof}

We are now ready to prove \cref{prop:posscalarssups}.

\begin{proof}[Proof of \cref{prop:posscalarssups}]
For existence, combine \cref{prop:commacatlimits,lem:normsq_preserve,lem:factornormsq}. As \(a = \normSq (X,x)\) for some object \((X,x)\) of \(\CommaCat{I}{U}\) by \cref{lem:epicrepresentative}, and the diagram
\[
\begin{tikzcd}[column sep=large]
\CommaCat{I}{U}
  \arrow[d, "{(X,x) \oplus
  \blank}" swap]
  \arrow[r, "\normSq"]
  \&
\PosScalars
  \arrow[d, "{a + \blank}"]
  \\
\CommaCat{I}{U}
  \arrow[r, "\normSq" swap]
  \&
\PosScalars
\end{tikzcd}
\]
commutes, preservation follows from \cref{prop:commacatlimits,lem:normsq_preserve}.
\end{proof}

\section{Recovering the real or complex numbers}
\label{sec:limits}

In operator algebra, a partial order is called \textit{monotone sequentially complete} (or \textit{monotone \(\sigma\)-complete}) if every bounded increasing sequence has a supremum. \textcite{demarr:partiallyorderedfields} showed that every partially ordered field that is monotone sequentially complete is isomorphic to~\(\Reals\). Whilst it is easy to define a partial order on~\(\SAScalars\) that is compatible with the field operations (see~\cref{lem:field_ordering} below), attempting to construct 
suprema of bounded increasing sequences in \(\SAScalars\) directly from nice category-theoretic assumptions like \cref{axiom:colimits:mono} seems futile. On the other hand, we already know that the partial order \(\leq\) on \(\PosScalars\) is monotone sequentially complete; perhaps there is a \citeauthor{demarr:partiallyorderedfields}-like theorem about partially ordered strict semifields that we might use instead? To answer this question, we need a better understanding of the properties of suprema and infima in partially ordered strict semifields.

First of all, due to the existence of multiplicative inverses, suprema and infima are always compatible with multiplication. 

\begin{proposition}[Compatibility with multiplication]
\label{lem:mult_extrema}
\label{lem:inf_mult_seq}
\label{lem:sup_mult_seq}
Let \(S\) be a partially ordered strict semifield. For all decreasing sequences \(a_n\) and \(b_n\) in \(S\),
\begin{enumerate}
    \item 
    \label{item:inf_mult}
    if \(\inf a_n\) and \(\inf b_n\) exist, then \[\inf a_nb_n = \inf a_n \inf b_n;\]
    
    \item
    \label{item:inf_quot}
    if \(\inf a_n b_n\) and \(\inf b_n\) exist and \(\inf b_n \neq 0\), then
    \[\inf a_n = \frac{\inf a_n b_n}{\inf b_n};\]

    \item
    \label{item:sup_inv_inf}
    if \(\inf b_n\) exists and \(\inf b_n \neq 0\), then \[\sup \frac{1}{b_n} = \frac{1}{\inf b_n}.\]
\end{enumerate}
Dually, for all increasing sequences \(a_n\) and \(b_n\) in \(S\),
\begin{enumerate}[resume]
    \item
    \label{item:sup_mult}
    if \(\sup a_n\) and \(\sup b_n\) exist then \[\sup a_nb_n = \sup a_n \sup b_n;\]
    
    \item
    \label{item:sup_quot}
    if \(\sup a_n b_n\) and \(\sup b_n\) exist, and \(b_1 \neq 0\), then \[\sup a_n = \frac{\sup a_n b_n}{\sup b_n};\]

    \item
    \label{item:inf_inv_sup}
    if \(\sup b_n\) exists and \(b_1 \neq 0\), then \[\inf \frac{1}{b_n} = \frac{1}{\sup b_n}.\]
\end{enumerate}
\end{proposition}

\begin{lemma}[Inversion is anti-monotonic]
\label{lem:inversion_antimonotonic}
In a partially ordered strict semifield, if \(a \leq b\) and \(a \neq 0\) then \(b \neq 0\) and \(\frac{1}{b} \leq \frac{1}{a}\).
\end{lemma}

\begin{proof}
If \(b = 0\) then \(0 \leq a \leq 0\) so \(a = 0\). Hence \(b \neq 0\) and \(\frac{1}{b} = \frac{a}{ab} \leq \frac{b}{ab} = \frac{1}{a}\).
\end{proof}

\begin{proof}[Proof of \cref{lem:mult_extrema}]
Let \(a_n\) and \(b_n\) be decreasing sequences in \(S\). If \(a_n\) is eventually zero then \cref{item:inf_mult} and  \cref{item:inf_quot} hold trivially, so assume that all \(a_k\) are non-zero.

For \cref{item:inf_mult}, suppose that \(\inf a_n\) and \(\inf b_n\) exist. For all~\(k\), we have \(\inf a_n \inf b_n \leq a_k b_k\). Suppose that \(c \leq a_k b_k\) for each \(k\). Then, for each~\(k\), as \(c/a_k \leq c/a_j \leq b_j\) for each \(j \geq k\), we have \(c/a_k \leq \inf b_n\), and thus \(c \leq a_k \inf b_n\). If \(\inf b_n = 0\), then \(c = 0 = \inf a_n \inf b_n\). If \(\inf b_n \neq 0\), then \(c /\inf b_n \leq a_k\) for each \(k\), so \(c /\inf b_n \leq \inf a_n\), and thus \(c \leq \inf a_n \inf b_n\).

For \cref{item:inf_quot}, suppose that \(\inf a_n b_n\) and \(\inf b_n\) exist, and that \(\inf b_n \neq 0\). For all~\(k\), we have \(\inf a_n b_n \leq a_j b_j \leq a_k b_j\) for each \(j \geq k\), so \(\inf a_nb_n \leq a_k \inf b_n\), and thus \(\inf a_n b_n /\inf b_n \leq a_k\). If \(c \leq a_k\) for each \(k\), then \(c \inf b_n \leq c b_k \leq a_k b_k\) for each~\(k\), so \(c \inf b_n \leq \inf a_n b_n\) and thus \(c \leq \inf a_n b_n / \inf b_n\).

For \cref{item:sup_inv_inf}, suppose that \(\inf b_n \neq 0\). For all \(k\), as \(\inf b_n \leq b_k\), also \(1/\inf b_n \geq 1/b_k\). If \(c \geq 1/b_k\) for all \(k\), then \(1/c \leq b_k\) for all \(k\), so \(1/c \leq \inf b_n\), and thus \(c \geq 1/\inf b_n\).

The dual statements about suprema may be proved similarly.
\end{proof}

\begin{proposition}
\label{lem:inf_sup_equivalence}
A partially ordered strict semifield is monotone sequentially complete if and only if it has infima of (bounded) decreasing sequences.
\end{proposition}

\begin{proof}
Suppose that it has infima of decreasing sequences. Let \(a_n\) be a bounded increasing sequence. If it is identically zero, then \(\sup a_n = 0\). Otherwise there is a \(j\) such that \(a_k \neq 0\) for all \(k \geq j\), so \(\sup a_n = \sup_n a_{n + j}\) exists by \cref{lem:mult_extrema}~\cref{item:sup_inv_inf}.

Suppose now that it has suprema of bounded increasing sequences. Let \(a_n\) be a decreasing sequence. If its only lower bound is \(0\), then \(\inf a_n = 0\). Otherwise it has a non-zero lower bound \(c\). Then all \(a_k\) are also non-zero. The sequence \(1/a_n\) has upper bound \(1/c\), so \(\inf a_n\) exists by \cref{lem:mult_extrema}~\cref{item:inf_inv_sup}.
\end{proof}

Due to the lack of additive inverses, compatibility of suprema and infima with addition is not guaranteed. Nevertheless, such compatibility is still quite a natural property of partially ordered strict semifields, holding, for example, in \(\PosReals\) and \(\TropReals\).

\begin{proposition}[Compatibility with addition]
\label{lem:add_extrema}
Let \(S\) be a partially ordered strict semifield. The following statements are equivalent:
\begin{enumerate}
    \item
    \label{item:sup_add_one}
    for all increasing sequences \(b_n\) in \(S\), if \(\sup b_n\) exists, then
    \[\sup(1 + b_n) = 1 + \sup b_n;\]
    
    \item
    \label{item:sup_add_seq}
    for all increasing sequences \(a_n\) and \(b_n\) in \(S\), if \(\sup a_n\) and \(\sup b_n\) exist, then
    \[\sup (a_n + b_n) = \sup a_n + \sup b_n;\]
    
    \item
    \label{item:inf_add_one}
    for all decreasing sequences \(b_n\) in \(S\), if \(\inf b_n\) exists and \(\inf b_n \neq 0\), then
    \[\inf(1 + b_n) = 1 + \inf b_n;\]
    
    \item
    \label{item:inf_add_seq}
    for all decreasing sequences \(a_n\) and \(b_n\) in \(S\), if \(\inf a_n\) and \(\inf b_n\) exist, and \(\inf a_n \neq 0\) and \(\inf b_n \neq 0\), then
    \[\inf (a_n + b_n) = \inf a_n + \inf b_n.\]
\end{enumerate}
\end{proposition}

\begin{proof}
Clearly \cref{item:sup_add_seq} implies \cref{item:sup_add_one}. For \cref{item:sup_add_one} implies \cref{item:sup_add_seq}, first observe that
\[
    a + \sup b_n = a\paren[\big]{1 + \sup \tfrac{b_n}{a}} = a \sup\paren[\big]{1 + \tfrac{b_n}{a}} = \sup(a + b_n)
\]
for all \(a \in S\). Hence
\[\sup\nolimits_m a_m + \sup\nolimits_n b_n = \sup\nolimits_m \paren[\big]{a_m + \sup\nolimits_n b_n} = \sup\nolimits_m \sup\nolimits_n(a_m + b_n).\]
For all \(k\), we have
\(\sup\nolimits_m \sup\nolimits_n(a_m + b_n) \geq \sup\nolimits_n(a_k + b_n) \geq a_k + b_k\). Suppose that \(c \geq a_\ell + b_\ell\) for all \(\ell\). For all \(j\) and \(k\), letting \(\ell = \max(j,k)\), we have \(c\geq a_\ell + b_\ell \geq a_j + b_k\). Hence \(c \geq \sup\nolimits_n(a_j + b_n)\) for all \(j\), and so
\(c \geq \sup\nolimits_m \sup\nolimits_n(a_m + b_n)\).

For \cref{item:sup_add_one} implies \cref{item:inf_add_one},
\[
1 + \inf b_n
= 1 + \frac{1}{\sup\frac{1}{b_n}}
= \frac{\sup \frac{1}{b_n} + 1}{\sup \frac{1}{b_n}}
= \frac{\sup \frac{1 + b_n}{b_n}}{\sup \frac{1}{b_n}}
= \frac{1}{\sup\frac{1}{1 + b_n}} = \inf(1 + b_n).
\]

Dually, \cref{item:inf_add_one} is equivalent to \cref{item:inf_add_seq}, and \cref{item:inf_add_one} implies \cref{item:sup_add_one}.
\end{proof}

A partially ordered strict semifield will be called \textit{suprema compatible} if it satisfies one of the equivalent conditions in \cref{lem:add_extrema}. 
Conditions \cref{item:inf_add_one,item:inf_add_seq} are still equivalent when \(\inf a_n\) and \(\inf b_n\) are allowed to be zero; a partially ordered strict semifield satisfying one of these stronger versions of conditions \cref{item:inf_add_one,item:inf_add_seq} will be called \textit{infima compatible}. By \cref{lem:add_extrema}, every infima-compatible partially ordered strict semifield is suprema compatible. The converse is not true.

\begin{example}
\label{ex:semifield}
The set \(\field{S} = \setb[\big]{\begin{psmallmatrix}x\\y\end{psmallmatrix} \in \Reals^2}{\text{\(x > 0\) and \(y > 0\), or \(x = y = 0\)}}\) is a partially ordered strict semifield with \(0 = \begin{psmallmatrix}0\\0\end{psmallmatrix}\), \(1 = \begin{psmallmatrix}1\\1\end{psmallmatrix}\), pointwise addition and multiplication, and \(\begin{psmallmatrix}x\\y\end{psmallmatrix} \leq \begin{psmallmatrix}u\\v\end{psmallmatrix}\) if and only if \(x \leq u\) and \(y \leq v\). It is monotone sequentially complete and suprema compatible, but not infima compatible. Indeed
\[
    \begin{psmallmatrix}1\\1\end{psmallmatrix} +  \inf \begin{psmallmatrix}1\\ 1/n\end{psmallmatrix}
    = \begin{psmallmatrix}1\\1\end{psmallmatrix} + \begin{psmallmatrix}0\\ 0\end{psmallmatrix}
    = \begin{psmallmatrix}1\\1\end{psmallmatrix}
    \neq
    \begin{psmallmatrix}2\\1\end{psmallmatrix}
    = \inf\begin{psmallmatrix}2\\1 + 1/n\end{psmallmatrix}
    = \inf \paren[\big]{\begin{psmallmatrix}1\\1\end{psmallmatrix} + \begin{psmallmatrix}1\\1/n\end{psmallmatrix}}.\]
\end{example}

This example shows that suprema compatibility is not enough to ensure that a partially ordered strict semifield that is monotone sequentially complete is isomorphic to~\(\PosReals\). On the other hand, infima compatibility, together with the inequality \(1 + 1 \neq 1\), is actually enough (see \cref{prop:pos_reals}). Unable to prove infima compatibility of \(\PosScalars\) directly from \cref{axiom:colimits:mono}, we need to use some additional property of \(\PosScalars\) to deal with the decreasing sequences that have infimum zero. 

The obvious candidate is the existence of an embedding into a field. Whilst \(\PosScalars\) embeds in a field, the semifield \(\field{S}\) from \cref{ex:semifield} does not. If it did, then
\[
  \paren[\big]{\begin{psmallmatrix} 2 \\ 1 \end{psmallmatrix} - \begin{psmallmatrix} 1 \\ 1 \end{psmallmatrix}}\paren[\big]{\begin{psmallmatrix} 2 \\ 1 \end{psmallmatrix} - \begin{psmallmatrix} 2 \\ 2 \end{psmallmatrix}} = \begin{psmallmatrix} 4 \\ 1 \end{psmallmatrix} - \begin{psmallmatrix} 2 \\ 1 \end{psmallmatrix} - \begin{psmallmatrix} 4 \\ 2 \end{psmallmatrix} + \begin{psmallmatrix} 2 \\ 2 \end{psmallmatrix} = \begin{psmallmatrix} 6 \\ 3 \end{psmallmatrix} - \begin{psmallmatrix} 6 \\ 3 \end{psmallmatrix} = \begin{psmallmatrix} 0 \\ 0 \end{psmallmatrix},
\]
and so \(\begin{psmallmatrix} 2 \\ 1 \end{psmallmatrix} = \begin{psmallmatrix} 1 \\ 1 \end{psmallmatrix}\) or \(\begin{psmallmatrix} 2 \\ 1 \end{psmallmatrix} = \begin{psmallmatrix} 2 \\ 2 \end{psmallmatrix}\). A similar trick yields the following proposition.

\begin{proposition}
\label{prop:inf_sum_geom}
Let \(S\) be a partially ordered strict semifield that is suprema compatible, monotone sequentially complete, and embeds in a field. For all \(a, u \in S\) with \(a \neq 0\) and \(u < 1\), we have \(\inf(a + u^n) = a\).
\end{proposition}

\begin{lemma}
In a partially ordered strict semifield, if \(a_n\) is a decreasing sequence and \(\inf a_n\) exists then \(\inf a_{2n} = \inf a_n\).
\end{lemma}

\begin{proof}
Firstly, \(\inf a_n \leq a_{2k}\) for each \(k\). Suppose that \(c \leq a_{2k}\) for each \(k\). Then, for each \(j\), either \(j = 2k\), in which case \(c \leq a_{2k} = a_j\), or \(j = 2k + 1\), in which case \(c \leq a_{2k + 2} \leq a_{2k + 1} = a_j\). Either way, it follows that \(c \leq \inf a_n\).
\end{proof}

\begin{proof}[Proof of \cref{prop:inf_sum_geom}]
As addition preserves non-zero infima,
\begin{align*}
    a + a^2 + \inf(a + u^n)^2
    &= \inf \paren[\big]{a + a^2 + (a + u^n)^2} = \inf \paren[\big]{2a(a + u^n) + (a + u^{2n})}
    \\&= 2a \inf (a + u^n) + \inf (a + u^{2n})
    = (2a + 1) \inf (a + u^n).
\end{align*}
Thus 
\(\paren[\big]{\inf(a + u^n) - a}\paren[\big]{\inf(a + u^n) - (a + 1)} = 0\) in the field.
If \(\inf(a + u^n) = a + 1\), then \(a + 1 \leq a + u \leq a + 1\), so \(a + u = a + 1\), and thus \(u = 1\), which is a contradiction. Hence \(\inf(a + u^n) = a\).
\end{proof}

The ideas above give rise to the following new characterisation of \(\Reals\) and \(\Comps\).

\begin{proposition}
\label{prop:comps_by_suprema}
Let \(C\) be an involutive field that has a partially ordered strict subsemifield \(P\) whose elements are all self-adjoint and include \(a^\dagger a\) for all \(a \in C\). If~\(P\) is suprema compatible and monotone sequentially complete, then there is an isomorphism of \(C\) with \(\Reals\) or \(\Comps\) that maps \(P\) onto \(\PosReals\).
\end{proposition}

Having already shown in \cref{prop:orderedsemifield,prop:posscalarssups} that \(\Scalars\) and \(\PosScalars\) satisfy the assumptions of \cref{prop:comps_by_suprema}, the following corollary is immediate.

\begin{corollary}
\label{thm:complexes}
There is an isomorphism of \(\Scalars\) with \(\Reals\) or \(\Comps\) that maps \(\PosScalars\) onto \(\PosReals\).
\end{corollary}

The remainder of this section is devoted to proving \cref{prop:comps_by_suprema}.

\begin{definition}[\textcite{demarr:partiallyorderedfields}]
A \textit{partially ordered field} is a field \(F\) equipped with a partial order \(\cleq\) satisfying the following axioms:
\begin{enumerate}
\item 
\label{axiom:pofield:add}
if $a \cleq b$ then $a+c \cleq b+c$,

\item
\label{axiom:pofield:mult}
if $0 \cleq a$ and $0 \cleq b$ then $0 \cleq ab$,

\item
\label{axiom:pofield:one}
$0 \cleq 1$,

\item
\label{axiom:pofield:inv}
if $0 \cleq a$ and $a \neq 0$ then $0 \cleq a^{-1}$, and

\item
\label{axiom:pofield:diff}
every $a \in F$ is of the form $a=b-c$ where $0 \cleq b$ and $0 \cleq c$.
\end{enumerate}
\end{definition}

\begin{lemma}
\label{lem:field_ordering}
Let \(R\) be a field with a partially ordered strict subsemifield \(P\) that contains all squares of \(R\). The binary relation \(\cleq\) on \(R\) defined by \(a \cleq b\) if \(b - a \in P\) is a partial order making \(R\) into a partially ordered field.
\end{lemma}

\begin{proof}
For reflexivity, \(a \cleq a\) because \(a - a = 0 \in P\). For antisymmetry, if \(a \cleq b\) and \(b \cleq a\), then \(a - b = 0\) because \(0 \leq a - b \leq (a - b) + (b - a) = 0\), so \(a = b\). For transitivity, if \(a \cleq b\) and \(b \cleq c\), then \(a \cleq c\) because \(c - a = (c - b) + (b - a) \in P\).

Axioms \cref{axiom:pofield:add,axiom:pofield:mult,axiom:pofield:one,axiom:pofield:inv} of a partially ordered field are straightforward to check. For example, if \(a \cleq b\), then \((b+c)-(a+c) = b - a \in P\), so \(a+c \cleq b+c\). For axiom~\cref{axiom:pofield:diff}, observe that \(a = (a+\frac{1}{2})^2 - (a^2+\frac{1}{4})\) where \(0 \cleq (a+\frac{1}{2})^2\) and \(0 \cleq a^2+\frac{1}{4}\) because \(P\) contains all squares of \(R\). 
\end{proof}

\begin{lemma}
\label{prop:inf_epsilon_pos_scalar}
Under the assumptions of \cref{prop:comps_by_suprema}, for each self-adjoint \(a \in C\), actually \(a \in P\) if and only if \(a + 2^{-k} \in P\) for each \(k\).
\end{lemma}

\begin{proof}
The \textit{only if} direction is trivial. For the \textit{if} direction, suppose that \(a + 2^{-k} \in P\) for each \(k\). Then \(\inf (a + 2^{-n})\) exists because
\[a + 2^{-k} = a + 2^{-(k + 1)} + 2^{-(k + 1)} \geq a + 2^{-(k + 1)}\]
for each \(k\). Either \(\inf(a + 2^{-n}) \neq 0\) or \(\inf(a + 2^{-n}) = 0\). In the former case,
\[
  \inf(a + 2^{-n}) + a^2 + \tfrac{1}{4}
  = \inf\paren[\big]{(a + \tfrac{1}{2})^2 + 2^{-n}}
  = \paren[\big]{a + \tfrac{1}{2}}^2
  = a^2 + a + \tfrac{1}{4},
\]
so \(a = \inf(a + 2^{-n}) \in P\). In the latter case,
\[0 = \inf(a + 2^{-n})^2 = \inf\paren[\big]{(a + 2^{-n})^2} = \inf\paren[\big]{a^2 + 2^{1-n}(a + 2^{-(n + 1)})}.\]
As \(0 \leq a^2 \leq a^2 + 2^{1-n}(a + 2^{-(n + 1)})\), it follows that
\[0 \leq a^2 \leq \inf\paren[\big]{a^2 + 2^{1-n}(a + 2^{-(n + 1)})} = 0,\]
so \(a^2 = 0\), and thus \(a = 0\). Hence \(a = \inf(a + 2^{-n}) \in P\) in this case as well.
\end{proof}

\begin{lemma}
\label{lem:sa_iso_reals}
Under the assumptions of \cref{prop:comps_by_suprema}, there is an isomorphism of the field \(R\) of self-adjoint elements of \(C\) with \(\Reals\) that maps \(P\) onto \(\PosReals\).
\end{lemma}

\begin{proof}
As \(P\) is a subsemifield of \(R\) containing all squares of \(R\), the binary relation \(\cleq\) defined in \cref{lem:field_ordering} makes \(R\) into a partially ordered field. By \citeauthor{demarr:partiallyorderedfields}'s theorem~\cite{demarr:partiallyorderedfields}, it suffices to show that \(\cleq\) is monotone sequentially complete.

Let \(a_1 \cleq a_2 \cleq \dots \cleq b\) in \(R\). Then $0 \leq a_1-a_1 \leq a_2-a_1 \leq a_3-a_1 \leq \cdots \leq b-a_1$, so \(\sup\nolimits_n(a_n - a_1)\) exists. 
We will show that \(a_n\) has supremum \(a = a_1 + \sup\nolimits_n(a_n - a_1)\) with respect to \(\cleq\). Firstly, for each \(k\), we have
\begin{align*}
&\nqquad\sup\nolimits_n\paren{a_{n + k - 1} - a_k} + {a_k}^2 + \tfrac{1}{4} + \paren{a_1 - \tfrac{1}{2}}^2\\
  &= \sup\nolimits_n\paren[\big]{a_{n + k - 1} + \paren{a_k - \tfrac{1}{2}}^2 + \paren{a_1 - \tfrac{1}{2}}^2}\\
  &= \sup\nolimits_n\paren{a_{n + k - 1} - a_1} + {a_1}^2 + \tfrac{1}{4} + \paren{a_k - \tfrac{1}{2}}^2\\
  &= \sup\nolimits_n\paren{a_n - a_1} + {a_1}^2 + \tfrac{1}{4} + \paren{a_k - \tfrac{1}{2}}^2\\
  &= a - a_1 + {a_1}^2 + \tfrac{1}{4} + \paren{a_k - \tfrac{1}{2}}^2
\end{align*}
so \(a = a_k + \sup_n(a_{n + k - 1} - a_k)\), and thus \(a_k \cleq a\).

Fix a natural number \(j\). The sequences \(a_n - a_1\) and \(b - a_n + 2^{-j}\) in \(P\) are respectively increasing and decreasing in \(n\). Now
\begin{align*}
b - a_1 + 2^{-j} 
    &= (a_k - a_1) + (b - a_k + 2^{-j}) \\
    &\leq \sup\nolimits_m(a_m - a_1) + b - a_k + 2^{-j}\\
    &= (a - a_1) + (b - a_k + 2^{-j})
\\
\shortintertext{for all \(k\), so}
b - a_1 + 2^{-j} 
    &\leq \inf\nolimits_n\paren[\big]{(a - a_1) + (b - a_n + 2^{-j})}\\
    &= a - a_1 + \inf\nolimits_n (b - a_n + 2^{-j})
\\\intertext{because \(\inf\nolimits_n(b - a_n + 2^{-j}) \geq 2^{-j} > 0\). Also}
b - a_1 + 2^{-j} 
    &= (a_k - a_1) + (b - a_k + 2^{-j})\\
    &\geq a_k - a_1 + \inf\nolimits_n (b - a_n + 2^{-j})\\
\shortintertext{for all \(k\), so}
b - a_1 + 2^{-j}
    &\geq \sup\nolimits_m \paren[\big]{a_m - a_1 + \inf\nolimits_n (b - a_n + 2^{-j})} \\
    &= \sup\nolimits_m (a_m - a_1) + \inf\nolimits_n (b - a_n + 2^{-j})\\
    &= a - a_1 + \inf\nolimits_n (b - a_n + 2^{-j}).\\
\intertext{Hence \(b - a_1 + 2^{-j} = a - a_1 + \inf\nolimits_n (b - a_n + 2^{-j})\), and so}
b - a + 2^{-j} &= \inf\nolimits_n (b - a_n + 2^{-j}) \in P.
\end{align*}
\Cref{prop:inf_epsilon_pos_scalar} now gives \(b - a \in P\), so \(a \cleq b\).
\end{proof}

\begin{lemma}
\label{lem:invol_reals_comps}
Let \(C\) be an involutive field, and \(R\) its subfield of self-adjoint elements. Each isomorphism \(\varphi \colon R \to \Reals\) such that \(\varphi(a^\dagger a) \geq 0\) for each \(a \in C\), uniquely extends to an isomorphism of \(C\) with \(\Reals\) or \(\Comps\).
\end{lemma}

\begin{proof}
The case \(C = R\) is trivial. Suppose that \(C \neq R\). Then there is a \(u \in C\) such that \(u \neq u^\dagger\). As \(\varphi \paren[\big]{(u - u^\dagger)^\dagger (u - u^\dagger)} > 0\), we may define \(r \in R \backslash \set{0}\) and \(i \in C\) by
\[
    r = \varphi^{-1} \sqrt{\varphi \paren[\big]{(u - u^\dagger)^\dagger (u - u^\dagger)}}
    \qquad\text{and}\qquad
    i = \frac{u - u^\dag}{r}.
\]
Then \(i^\dag = -i\) because \(r\) is self-adjoint, and \(i^2 = -1\) because \(r^2 = (u - u^\dagger)^\dagger (u - u^\dagger)\). It follows also that \(i \neq 0\) and so \(i^\dag \neq i\).

We now show that $\set{1,i}$ is a basis for $C$ as a vector space over $R$. It is a spanning set because every $a \in C$ satisfies the equation
\[
a = \frac{a + a^\dag}{2} + \frac{a - a^\dag}{2i}i,
\]
where $(a + a^\dag)/2$ and $(a - a^\dag)/2i$ are self-adjoint. For linear independence, let $a,b \in R$ and suppose that $a + bi = 0$. If $b \neq 0$ then $i = -a/b$ would be self-adjoint, which is a contradiction. Thus $b = 0$ and $a = -bi = 0$.

Define \(\psi \colon C \to \Comps\) by \(\psi(a + bi) = \varphi(a) + \varphi(b)i\) for all \(a, b \in R\). Using the equations \(i^2 = -1\) and \(i^\dagger = -i\), and the fact that \(\varphi\) is an isomorphism of fields, it is easy to check that \(\psi\) is an isomorphism of involutive fields that extends \(\varphi\).
\end{proof}

\begin{proof}[Proof of \cref{prop:comps_by_suprema}] Combine \cref{lem:sa_iso_reals,lem:invol_reals_comps}.
\end{proof}

\section{Finite dimensionality}\label{sec:finite}

Our goal now is to prove the following abstract version of \cref{prop:con_dagger_finite}, from which it easily follows that the inner-product space associated to each dagger-finite object is finite dimensional. Here, the term dagger finite is unambiguous because the dagger monomorphisms in \(\C\) all come from \(\D\)~\cite[Lemma~14]{heunenkornellvanderschaaf:con}.

\begin{proposition}
\label{prop:finite-dagger-biproduct}
  Every dagger-finite object is dagger isomorphic to one of the form
  \[I^{\oplus n} = \underbrace{I \oplus I \oplus \dots \oplus I}_{n \text{ times}}.\]
\end{proposition}

To prove this proposition, we will use the following two lemmas.

\begin{lemma}\label{lem:dagger-sub-dagger-finite}
  Let \(m \colon A \to X\) be a dagger monomorphism. If \(X\) is dagger finite, then \(A\) is also dagger finite.
\end{lemma}
\begin{proof}
Let \(m^\perp \colon A^\perp \to X\) be a kernel of~\(m^\dagger\). Then \(X\) is dagger isomorphic to \(A \oplus A^\perp\). Suppose that \(X\), and thus \(A \oplus A^\perp\), is dagger finite. Let \(f \colon A \to A\) be a dagger monomorphism. Then \(f \oplus 1 \colon A \oplus A^\perp \to A \oplus A^\perp\) is also dagger monic, and so a dagger isomorphism. Hence \(f\) is also a dagger isomorphism.
\end{proof}

\begin{lemma}\label{lem:nonzerovector}
  Every non-zero object \(X\) admits a dagger monomorphism from \(I\).
\end{lemma}
\begin{proof}
  Let \(X\) be a non-zero object. Then the morphisms \(0 \colon X \to X\) and \(1 \colon X \to X\) are distinct. As \(I\) is a separator, there is a morphism \(x \colon I \to X\) in \(\C\) such that \( 1x \neq 0x\), that is, such that \(x\) is non-zero. As \(\C\) has a zero object and dagger equalisers, the positive scalar \(x^\dagger x\) is then also non-zero~\cite[Lemma~II.5]{vicary:complex}, and so, via the isomorphism \(\PosScalars \cong \PosReals\), has a non-zero positive square root \((x^\dagger x)^{1/2}\). Let \(m = x (x^\dagger x)^{-1/2}\). Then \(m\colon I \to X\) is dagger monic because
  \(m^\dagger m = (x^\dagger x)^{-1/2} x^\dagger x (x^\dagger x)^{-1/2} = 1\).
\end{proof}

\begin{proof}[Proof of \cref{prop:finite-dagger-biproduct}]
  Let \(X\) be a dagger-finite object. An \textit{orthonormal system} of elements of \(X\) is a set of dagger monomorphisms \(\daggerMono{I}{X}\) that are pairwise orthogonal. As \(\D\) is locally small, and the dagger monomorphisms into \(X\) are in bijection with the dagger idempotents on~\(X\), the class of orthonormal systems of elements of \(X\) is itself a set. With respect to subset inclusion, the union of a chain of orthonormal systems of \(X\) is again an orthonormal system, so, by Zorn's lemma, there is an orthonormal system \(\mathcal{S}\) of \(X\) that is maximal.

  Assume that \(\mathcal{S}\) is infinite. Then it has a countable subset \(\set{x_k}_{k = 1}^\infty\). The dagger monomorphisms \(\begin{bsmallmatrix} x_1 & x_2 & \dots & x_k \end{bsmallmatrix} \colon I^{\oplus k} \to X\) form a cocone on the sequential diagram
  \[
  \begin{tikzcd}[sep=large, cramped]
  I^{\oplus 1}
      \arrow[r, "i_{1}"]
      \&[-1em]
  I^{\oplus 2}
      \arrow[r, "i_{1,2}"]
      \&
  I^{\oplus 3}
      \arrow[r, "i_{1,2,3}"]
      \&[1em]
  \cdots
  \end{tikzcd}
  \]
  in \(\D\), so, by \cref{axiom:colimits:mono}, this diagram has a colimit \(j_k \colon I^{\oplus k} \to I^{\oplus \Nats}\). This colimit is also a colimit in the wide subcategory of \(\D\) of dagger monomorphisms~\cite[Lemma~20]{heunenkornellvanderschaaf:con}.\footnote{The arXiv version of the article proves Lemma 20 in more detail than the published version.} There is thus a unique dagger monomorphism \(x \colon I^{\oplus \Nats} \to X\) such that \(x_k = xj_k\) for each natural number \(k\). By \cref{lem:dagger-sub-dagger-finite}, the object \(I^{\oplus \Nats}\) is also dagger finite. There is also a unique dagger monomorphism \(s \colon I^{\oplus \Nats} \to I^{\oplus \Nats}\) such that the diagram
  \[
  \begin{tikzcd}[sep=large]
  I^{\oplus 1}
      \arrow[r, "i_1" swap]
      \arrow[d, "i_2" swap]
      \arrow[rrr, "j_1", bend left=20, shorten >=1ex, shift left]
      \&
  I^{\oplus 2}
      \arrow[r, "i_{1,2}" swap]
      \arrow[d, "i_{2,3}" swap]
      \arrow[rr, "j_2", bend left=10, shift left]
      \&[1em]
  \cdots
      \&
  I^{\oplus \Nats}
      \arrow[d, "s"]
  \\
  I^{\oplus 2}
      \arrow[r, "i_{1,2}"]
      \arrow[rrr, "j_2" swap, bend right=20, shorten >=2ex, shift right]
      \&
  I^{\oplus 3}
      \arrow[r, "i_{1,2,3}"]
      \arrow[rr, "j_3" swap, bend right=10, shift right]
      \&[1em]
  \cdots
      \&
  I^{\oplus \Nats}
  \end{tikzcd}
  \]
  commutes. The morphism \(s\) is actually a dagger isomorphism by \cref{axiom:finite}. Now
\[
\begin{tikzcd}[column sep={scriptsize}, row sep={large}]
    \&
    I^{\oplus \Nats}
        \arrow[dr, "s"]
    \&
    \&
\\
    I^{\oplus k}
        \arrow[ur, "j_k"]
        \arrow[dr, "i_{2,\dots,k + 1}"]
        \arrow[rrrr, "0" swap, out=-80, in=-100, looseness=1.2]
    \&
    \&
    I^{\oplus \Nats}
        \arrow[rr, "{j_1}^\dagger"]
        \arrow[dr, "{j_{k + 1}}^\dagger"]
    \&
    \&
    I^{\oplus 1}
\\
    \&
    I^{\oplus(k + 1)}
        \arrow[rr, "1" swap]
        \arrow[ur, "j_{k + 1}" swap]
    \&
    \&
    I^{\oplus(k + 1)}
        \arrow[ur, "{i_1}^\dagger"]
\end{tikzcd}
\]
commutes
  % \[{j_1}^\dagger s j_k = {j_1}^\dagger j_{k + 1} i_{2, 3, \dots, k + 1} = {i_1}^\dagger {j_{k + 1}}^\dagger j_{k + 1} i_{2, 3, \dots, k + 1} = {i_1}^\dagger i_{2, 3, \dots, k + 1} = 0\]
  for each natural number \(k\). As the colimit cocone \(j_k\) is jointly epic in \(\D\), it follows that \({j_1}^\dagger s = 0\), and so \({j_1}^\dagger j_1 = {j_1}^\dagger ss^\dagger j_1 = 0 \neq 1\), which is a contradiction.

  Hence \(\mathcal{S}\) is finite, and so of the form \(\set{x_k}_{k = 1}^n\). If the orthogonal complement of \(\begin{bsmallmatrix}x_1 & x_2 & \cdots & x_n\end{bsmallmatrix} \colon I^{\oplus n} \to X\) were non-zero, then, using \cref{lem:nonzerovector}, we could obtain a dagger monomorphism \(I \to X\) that is orthogonal to all \(x_k\), contradicting maximality of \(\mathcal{S}\). Hence \(\begin{bsmallmatrix}x_1 & x_2 & \cdots & x_n\end{bsmallmatrix}\) is a dagger isomorphism.
\end{proof}

\begin{proposition}\label{prop:fhilb}
The dagger rig categories \(\C\) and \(\FHilb_{\Scalars}\) are equivalent.
\end{proposition}

\begin{proof}
  As the monoidal structure $(\oplus,O)$ on $\C$ is a choice of zero object and binary dagger biproducts, and these are preserved by equivalences of dagger categories, it suffices to construct an equivalence of dagger monoidal categories with respect to \(\otimes\).

  Let \(X\) be an object of \(\C\). The equation
  \(\innerProd{x_1}{x_2} = {x_1}^\dagger x_2\)
  defines a \(\Scalars\)-valued inner product on \(\C(I, X)\). The object \(X\) is a dagger biproduct in \(\C\) of finitely many copies of \(I\) by \cref{axiom:finite,prop:finite-dagger-biproduct}. The biproduct injections form an orthonormal basis for \(\C(I, X)\), so this inner-product space is actually finite dimensional.

  For each morphism \(f \colon X \to Y\) in \(\C\), the function \(\C(I, f)\) has adjoint \(\C(I, f^\dagger)\) because
  \(\innerProd{y}{fx} = y^\dagger fx = y^\dagger f^{\dagger\dagger} x = (f^\dagger y)^\dagger x = \innerProd{f^\dagger y}{x}\), and so is also \(\Scalars\)-linear.

  We now know that the functor \(\C(I,\blank) \colon \C \to \Set\) corestricts along the forgetful functor \(\FHilb_{\Scalars} \to \Set\), and that the corestriction is a dagger functor. It is actually dagger strong monoidal~\cite[Lemma~9]{heunenkornell:hilb}. Each Hilbert space of dimension \(n\) is dagger isomorphic to the Hilbert space $\C(I, I^{\oplus n})$, so it is dagger essentially surjective. It is full by the matrix calculus~\cite[Corollary~2.27]{heunenvicary:cqm} for $\C$ and $\FHilb_{\Scalars}$, and faithful because \(I\) is a separator.
\end{proof}

\section{Characterising \texorpdfstring{\(\FCon\)}{FCon} and \texorpdfstring{\(\Con\)}{Con}}
\label{sec:con}

In this brief final section, we finish proving our characterisation of \(\FCon\) and also sketch a new proof of the characterisation~\cite[Theorem~26]{heunenkornellvanderschaaf:con} of $\Con$ that bypasses Solèr's theorem.

The following statement of \cref{thm:main} is a rewording of the one in \cref{sec:axioms}.

\begingroup
\def\thetheorem{\ref*{thm:main}}
\begin{theorem}
  The dagger rig category \(\D\) is equivalent to the dagger rig category $\FCon$ of finite-dimensional Hilbert spaces and linear contractions.
\end{theorem}
\addtocounter{theorem}{-1}
\endgroup

\begin{proof}
  By \cref{prop:fhilb}, the dagger rig category $\D$ is equivalent to a wide dagger rig subcategory of $\FHilb$. To show that \(\D\) is equivalent to \(\FCon\), we must show that the morphisms of each are included in the morphisms of the other.

  By Sz.-Nagy's dilation theorem~\cite[Theorem~I.4.1]{nagy:harmonic}, each morphism \(f \colon X \to Y\) in \(\FCon\) has a factorisation \(f = em\) where \(m\) and \(e^\dagger\) are dagger monic. Both \(m\) and \(e\) come from \(\D\)~\cite[Lemma~14]{heunenkornellvanderschaaf:con}. As \(f\) is their composite, it also comes from \(\D\).
  
  Now let $f \colon X \to Y$ be a morphism in $\D$. To show that $f$ is in $\FCon$, it suffices to show that $\abs{\innerProd{fx}{y}} \leq 1$ for all $x \colon I \to X$ and $y \colon I \to Y$ in \(\C\) with norm \(1\). Such \(x\) and \(y\) are dagger monic and so actually come from $\D$. Hence $\innerProd{fx}{y} = x^\dag f^\dag y$ is a scalar in $\D$. But every scalar $z \colon I \to I$ in $\D$ satisfies $\abs{z}^2 = z^\dag z \leq 1$ by \cref{lem:monosinD}.
\end{proof}

The original proof of the characterisation~\cite[Theorem~26]{heunenkornellvanderschaaf:con} of \(\Con\) is spread across the articles~\cite{heunenkornellvanderschaaf:con} and~\cite{heunenkornell:hilb}, with passage from the former article to the latter one occurring in~\cite[Proposition~21]{heunenkornellvanderschaaf:con}. The only place in this proof where Solèr's theorem is used is in~\cite[Proposition~5]{heunenkornell:hilb}, which, in the language of the present article, says
\begin{enumerate}
    \item
    the involutive field \(\Scalars\) is isomorphic to \(\Reals\) or \(\Comps\), and
    \item for each object \(X\) of \(\C\), the \(\Scalars\)-inner product space \(\C(I, X)\) is actually a Hilbert space. 
\end{enumerate}
These facts can now be proved without Solèr's theorem. Fact~(1) is \cref{thm:complexes}. Fact~(2) follows from~\cite[Theorem~3.1]{gudder:innerproduct}, \cite[Lemma~4]{heunenkornell:hilb} and fact (1).

\section*{Declarations}

This version of the article has been accepted for publication, after peer review but is not the Version of Record and does not reflect post-acceptance improvements, or any corrections. The Version of Record is available online at the following URL:
\begin{center}
\url{http://dx.doi.org/10.1007/s10485-025-09803-5}.
\end{center}

\subsection*{Acknowledgements}

We are grateful to Andre Kornell for many discussions about this project, and particularly for the idea behind the proof of \cref{prop:finite-dagger-biproduct}.

{
\linespread{1.0}
\printbibliography
}

\appendix

\section{Completeness via colimits of epimorphisms}
\label{sec:colimits:epi}

In this appendix, replacing \cref{axiom:colimits:mono} with the following alternative completeness axiom, we show that the scalars are again the real or complex numbers. Our proof combines several ideas from across the literature~\cite{demarr:partiallyorderedfields,prestel:soler,fritz:semifield}. It is unclear whether this alternative axiom is strong enough to prove \cref{thm:main}.

\begin{alternativeaxiom}{axiom:colimits:mono}\label{axiom:colimits:epi}
  Every sequential diagram of epimorphisms has a colimit, and, for each natural transformation of such diagrams whose components are dagger monic, the induced morphism between the colimits is also dagger monic.
\end{alternativeaxiom}

To illustrate this, let \(m_n \colon (X_n, f_n) \to (Y_n, g_n)\) be such a natural transformation, let $c_n \colon X_n \to \colim X_n$ and $d_n \colon Y_n \to \colim Y_n$ be colimit cocones on \((X_n, f_n)\) and \((Y_n, g_n)\), and let \(m_\infty \colon \colim X_n \to \colim Y_n\) be the unique morphism such that \(m_\infty c_n = d_n m_n\). \Cref{axiom:colimits:epi} says that $m_\infty$ is dagger monic if each $m_n$ is dagger monic.

\[
  \begin{tikzcd}[sep=large]
  X_1
      \arrow[r, "f_1" swap, epi]
      \arrow[d, "m_1" swap, dagger mono]
      \arrow[rrr, "c_1", bend left=20, shorten >=1ex, shift left]
      \&
  X_2
      \arrow[r, "f_2" swap, epi]
      \arrow[d, "m_2" swap, dagger mono]
      \arrow[rr, "c_2", bend left=10, shift left]
      \&[1em]
  \cdots
      \&
 \colim X_n
      \arrow[d, "m_\infty", dagger mono]
  \\
  Y_1
      \arrow[r, "g_1", epi]
      \arrow[rrr, "d_1" swap, bend right=20, shorten >=2ex, shift right]
      \&
  Y_2
      \arrow[r, "g_2", epi]
      \arrow[rr, "d_2" swap, bend right=10, shift right]
      \&[1em]
  \cdots
      \&
  \colim Y_n
  \end{tikzcd}
  \]

\begin{proposition}
  The category $\FCon$ satisfies \cref{axiom:colimits:epi}.
\end{proposition}
% \cref{axiom:colimits:epi} actually also holds for \(\Con\), and even for directed colimits rather than sequential ones. We do not prove these claims as the proofs are more complicated than the one above and the claims are not used.
The construction~\cite{heunenkornellvanderschaaf:con} of directed colimits in $\Con$ also works for sequential colimits of epimorphisms in $\FCon$. We describe a simpler construction possible in this case.
\begin{proof}
  Let \((X_n, f_n)\) be a sequential diagram of epimorphisms in \(\FCon\). For each \(x_1 \in X_1\), as the sequence
  \(\norm{x_1}, \norm{f_1 x_1}, \norm{f_2f_1 x_1}, \dots\) is decreasing, its limit exists and is equal to its infimum. The set
  \[N = \setb*{x_1 \in X_1}{\lim_{n \to \infty} \norm{f_n \dots f_2 f_1 x_1} = 0}\]
  is a vector subspace of \(X_1\). Let \(\colim(X_n)\) be the vector space \(X_1/N\) equipped with the inner product defined by the limit
  \[\innerProd{x_1 + N}{x_1' + N} = \lim_{n \to \infty} \innerProd{f_n \dots f_2f_1 x_1}{f_n \dots f_2f_1 x_1'},\]
  which exists by the polarisation identity. It is finite dimensional by construction.

  To define \(c_k \colon X_k \to \colim(X_n)\), let \(x_k \in X_k\). As \(f_{k-1} \dots f_2 f_1\) is an epimorphism of finite-dimensional Hilbert spaces, it is surjective, so there is an \(x_1 \in X_1\) such that \(x_k = f_{k-1} \dots f_2f_1 x_1\). Let \(c_k x_k = x_1 + N\). As \(f_{k-1} \dots f_2f_1 x_1 = f_{k-1} \dots f_2 f_1 x_1'\) implies \(x_1' - x_1 \in N\), the map \(c_k\) is well defined. The maps \(c_k \colon X_k \to \colim(X_n)\) form a colimit cocone on \((X_n, f_n)\).

  Let \((Y_n, g_n)\) be another sequential diagram of epimorphisms in \(\FCon\). Consider a natural transformation \(m_n \colon (X_n, f_n) \to (Y_n, g_n)\) with dagger monic components. Let \(m_\infty\) be the morphism \(\colim X_n \to \colim Y_n\) that it induces. 
  For each \(x_1 \in X_1\), 
  \begin{multline*}
      \norm{m_\infty c_1 x_1} = \norm{d_1m_1 x_1} = \lim_{k \to \infty} \norm{g_k \dots g_2 g_1 m_1 x_1} \\
      = \lim_{k \to \infty} \norm{m_{k + 1} f_n \dots f_2 f_1 x_1} = \lim_{k \to \infty} \norm{f_k \dots f_2 f_1 x_1} = \norm{c_1 x_1}
  \end{multline*}
  because each \(m_k\) is an isometry. As each element of \(\colim X_n\) is of the form \(c_1 x_1\), it follows that \(m_\infty\) is also an isometry.
\end{proof}

\subsection{Characterising the positive reals}
We will use the following proposition to deduce that \(\PosScalars\) is isomorphic to \(\PosReals\).
Recall that a partially ordered strict semifield is \emph{monotone sequentially complete} if it has infima of decreasing sequences, and is \textit{infima compatible} if every decreasing sequence \(b_n\) satisfies \(\inf (1 + b_n) = 1 + \inf b_n\).

\begin{proposition}\label{prop:pos_reals}
  A partially ordered strict semifield is isomorphic to \(\PosReals\) if and only if it is monotone sequentially complete, infima compatible, and \(1 + 1 \neq 1\).
\end{proposition}

This result is a new variant of~\cite[Theorem~4.5]{fritz:semifield}. The proof, spread across several lemmas below, is inspired by \citeauthor{demarr:partiallyorderedfields}'s proof of a similar result about partially ordered fields~\cite{demarr:partiallyorderedfields}. The main idea is that \(u < 1\) exactly when the geometric series \(1 + u + u^2 + \dots\) converges.

\begin{lemma}
\label{lem:total}
If a partially ordered strict semifield is monotone sequentially complete and infima compatible, then it is totally ordered.
\end{lemma}

\begin{proof}
  We must show that \(a \geq b\) or \(a \leq b\). If \(b = 0\), then \(a \geq 0 = b\). Suppose that \(b \neq 0\), and let \(u = ab^{-1}\). It suffices to show that \(u \geq 1\) or \(u \leq 1\).

  For each \(n\), let ${s_n} = u^{n} + \dots + u + 1$. Then \(s_n\) is increasing, so \(\frac{1}{s_n}\) is decreasing, and thus \(\inf \frac{1}{s_n}\) exists. Observe that
  \begin{equation}\tag{$*$}\label{eq:geom_series}
      u + \frac{1}{s_n} = \frac{us_n + 1}{s_n} = \frac{s_{n + 1}}{s_n} = \frac{u^{n + 1} + s_n}{s_n} = \frac{u^{n + 1}}{s_n} + 1.
  \end{equation}
  In particular, \(\frac{1}{s_n} = \paren[\big]{u + \frac{1}{s_n}}\tfrac{1}{s_{n + 1}}\), so, by \cref{lem:inf_mult_seq},
  \[\inf \frac{1}{s_n} = \inf\paren[\Big]{u + \frac{1}{s_n}}\inf \frac{1}{s_{n + 1}} = \paren[\Big]{u + \inf \frac{1}{s_n}}\inf \frac{1}{s_n}.\]
  If \(\inf \frac{1}{s_n} \neq 0\), then \(u + \inf \frac{1}{s_n} = 1\), and so \(u \leq 1\). Otherwise \(\inf \frac{1}{s_n} = 0\), so, by \cref{eq:geom_series},
  \[u = u + \inf \frac{1}{s_n} = \inf\paren[\Big]{u + \frac{1}{s_n}}= \inf\paren[\Big]{\frac{u^{n + 1}}{s_n} + 1} \geq 1.\qedhere\]
\end{proof}

A partially ordered strict semifield is called \emph{multiplicatively Archimedean} if \(a \leq 1\) whenever the set \(\set{1, a, a^2, \dots}\) has an upper bound~\cite[Definition~4.1]{fritz:semifield}.

\begin{lemma}\label{lem:archimedean}
  If a partially ordered strict semifield is monotone sequentially complete and infima compatible, then it is multiplicatively Archimedean.
\end{lemma}
\begin{proof}
  By \cref{lem:total}, the order is total. Suppose, for the contrapositive, that \(a > 1\). The sequence \(a^{-n}\) is decreasing, so \(\inf a^{-n}\) exists. Now \(\inf a^{-n} = 0\) because
  \[\inf a^{-n} = aa^{-1} \inf a^{-n} = a \inf a^{-(n + 1)} = a \inf a^{-n}\]
  and \(a \neq 1\). It follows that \(\set{1, a, a^2, \dots}\) has no upper bound.
\end{proof}

Before proving \Cref{prop:pos_reals}, observe that the isomorphisms of totally ordered semifields are those morphisms that are bijective. The monotonicity of the inverse function follows from totality of the order.

\begin{proof}[{Proof of \cref{prop:pos_reals}}]
Let \(S\) be a partially ordered strict semifield that is infima compatible and monotone sequentially complete. By \Cref{lem:total,lem:archimedean}, it embeds in~\(\PosReals\) or \(\TropReals\)~\cite[Theorem~4.2]{fritz:semifield}. As \(1 + 1 \neq 1\), the latter is impossible, so there is an embedding \(\varphi \colon S \hookrightarrow \PosReals\). We will show that \(\varphi\) is surjective, and so an isomorphism.

Firstly, every embedding into \(\PosReals\) is necessarily surjective on \(\PosRats\). Let \(r \in \PosReals\). Then \(\sup a_n = r = \inf b_n\) for some \(a_1 \leq a_2 \leq \dots \leq r \leq \cdots \leq b_2 \leq b_1\) in~\(\PosRats\). As \(S\) is monotone sequentially complete, the infimum \(\inf b_n\) also exists in \(S\). We will show that \(r = \varphi(\inf b_n)\). For all \(k\), we have \(\inf b_n \leq b_k\), so \(\varphi(\inf b_n) \leq \varphi(b_k) = b_k\). Hence \(\varphi(\inf b_n) \leq \inf b_n = r\). Also, for all \(k\), as \(a_k \leq b_j\) for all \(j\), also \(a_k \leq \inf b_n\), and thus \(a_k = \varphi(a_k) \leq \varphi(\inf b_n)\). Hence \(r = \sup a_n \leq \varphi(\inf b_n)\), so actually \(r = \varphi(\inf b_n)\).
\end{proof}

\subsection{The positive scalars are the positive reals}
To show that \(\PosScalars\) is isomorphic to~\(\PosReals\), we now check the hypotheses of \cref{prop:pos_reals}. The first part of \cref{axiom:colimits:epi} corresponds to monotone sequential completeness while the second part corresponds to infima compatibility.

\begin{lemma}\label{lem:comma_colim}
  The category \(\CommaCat{I}{U}\) has colimits of sequential diagrams of epimorphisms and, for each object \((X,x)\) of \(\CommaCat{I}{U}\), the endofunctor \((X,x) \oplus \blank\) preserves them.
\end{lemma}

\begin{proof}
The proofs of \cref{lem:biproductspreservecolimits} and \cref{prop:commacatlimits} work \textit{mutatis mutandis}, with \cref{prop:picreateslimits} replaced by the following reasoning. As the functor \(\Pi \colon \CommaCat{I}{U} \to \D\) creates connected colimits~\cite[Theorem~3]{burstall:computational} (see also \cite[Proposition~3.3.8]{riehl:categorycontext}), it creates epimorphisms, and thus also colimits of sequential diagrams of epimorphisms. 
\end{proof}

\begin{lemma}
\label{lem:normSq_preserve}
The functor \(\normSq \colon \CommaCat{I}{U} \to \PosScalars\) preserves colimits of sequential diagrams of epimorphisms.
\end{lemma}

\begin{proof}
A \textit{coimage} of a morphism \(f\) is a terminal epimorphism \(e\) through which \(f\) factors. This means that \(f = ge\) for some morphism \(g\), and, if \(f = g'e'\) and \(e'\) is epic, then \(e = h'e'\) for some morphism \(h'\) that is uniquely determined by the fact that \(e'\) is epic. Recall that every morphism \(f\) of \(\C\) has a factorisation \(f = me\) where \(m\) is dagger monic and \(e\) is epic. The epic part of such a factorisation of \(f\) is a coimage of~\(f\). Indeed, if \(f = g'e'\) and \(e'\) is epic, then \(e = m^\dagger me = m^\dagger f = m^\dagger g'e'\) so \(h' = m^\dagger g'\) satisfies \(e = h'e'\). For each object \((X,x)\) of \(\CommaCat{I}{U}\), choose an epimorphism \(e \colon \epi{I}{E}\) and a dagger monomorphism \(m \colon \daggerMono{E}{X}\) such that \(x = me\), and let \(\Coim (X,x) = (E,e)\) and \(\epsilon_{(X,x)} = m \colon J\Coim (X,x) \to (X,x)\), noting that \(m\) is actually from \(\D\)~\cite[Lemma~14]{heunenkornellvanderschaaf:con}. By the universal property of coimages, this data uniquely extends to a right adjoint \(\Coim\) of \(J\) with counit \(\epsilon\)~\cite[\S IV.1 Theorem~2]{maclane:1968:categories}. 

Now \(\normSq = \normSq \circ J \circ \Coim\) because the components of \(\epsilon\) are dagger monic. Also \(\normSq \circ J\) preserves all colimits by \cref{lem:normSq:factors}. Hence it suffices to show that colimits of sequential diagrams of epimorphisms are preserved by \(\Coim\). As \(J\) is a full subcategory inclusion with a right adjoint, it is actually a coreflective subcategory inclusion, so it creates all colimits that \(\CommaCat{I}{U}\) admits~\cite[Proposition~4.5.15]{riehl:categorycontext}. Thus to show that colimits of sequential diagrams of epimorphisms are preserved by \(\Coim\), it suffices to show that they are preserved by \(J \circ \Coim\).

Let \(\paren[\big]{(X_n, x_n), f_n}\) be a sequential diagram of epimorphisms in \(\CommaCat{I}{U}\), and let \(c_n \colon (X_n, x_n) \to (X_\infty, x_\infty)\) be a colimit cocone on this diagram. For each~\(k\), let \((E_k, e_k) = J\Coim(X_k, x_k)\), \(g_k = J\Coim f_k\), and \(m_k = \epsilon_{(X_k, x_k)}\). As \(\paren[\big]{(E_n, e_n), g_n}\) is also a sequential diagram of epimorphisms in \(\CommaCat{I}{U}\), it has a colimit cocone \(d_n \colon (E_n, e_n) \to (E_\infty, e_\infty)\) by \cref{lem:comma_colim}. Let \(m_\infty \colon (E_\infty, e_\infty) \to (X_\infty, x_\infty)\) be the unique morphism such that \(m_\infty d_k = c_k m_k\) for all \(k\).
\[
  \begin{tikzcd}
  (E_1, e_1)
      \arrow[r, "g_1" swap, epi]
      \arrow[d, "m_1" swap]
      \arrow[rrr, "d_1", bend left=20, shorten >=1ex, shift left]
      \&
  (E_2, e_2)
      \arrow[r, "g_2" swap, epi]
      \arrow[d, "m_2" swap]
      \arrow[rr, "d_2", bend left=10, shift left]
      \&[1em]
  \cdots
      \&
 (E_\infty, e_\infty)
      \arrow[d, "m_\infty"]
  \\
  (X_1, x_1)
      \arrow[r, "f_1", epi]
      \arrow[rrr, "c_1" swap, bend right=20, shorten >=2ex, shift right]
      \&
  (X_2, x_2)
      \arrow[r, "f_2", epi]
      \arrow[rr, "c_2" swap, bend right=10, shift right]
      \&[1em]
  \cdots
      \&
  (X_\infty, x_\infty)
  \end{tikzcd}
\]
As \(\Pi \colon \CommaCat{I}{U} \to \D\) creates connected limits, \(c_n \colon X_n \to X_\infty\) and \(d_n \colon E_n \to E_\infty\) are colimit cocones on the underlying sequential diagrams \((X_n, f_n)\) and \((E_n, g_n)\) in \(\D\), respectively. Thus \(m_\infty \colon E_\infty \to X_\infty\) is dagger monic by \cref{axiom:colimits:epi}.

We now show that the morphism \(e_\infty \colon I \to E_\infty\) in \(\C\) is epic. Let \(s, t \colon E_\infty \to A\) in \(\C\) and suppose that \(se_\infty = te_\infty\). Then, for all \(k\), as \(sd_ke_k = se_\infty = te_\infty = td_ke_k\) and \(e_k\) is epic, actually \(sd_k = td_k\). Now the cocone \(d_n \colon E_n \to E_\infty\) in \(\D\), being a colimit cocone, is jointly epic in \(\D\). Similarly to \cref{lem:monos}, the functor \(U \colon \D \to \C\) preserves jointly epic wide cospans. Hence \(d_n\) is also jointly epic in \(\C\). Thus \(s = t\).

As \(x_\infty = m_\infty e_\infty\) is an (epic, dagger monic) factorisation of \(x_\infty\), there is a unique isomorphism \(u \colon (E_\infty, e_\infty) \to \Coim (X_\infty, x_\infty)\). Also, for all \(k\), the morphism in \(\D\) underlying \(\Coim c_k\) is equal to \(u d_k\), by the universal property of coimages. Hence \(J \Coim c_n\) is another colimit cocone on \(\paren[\big]{(E_n, e_n), g_n}\).
\end{proof}

\begin{proposition}
\label{prop:pos_reals_concrete}
  The partially ordered strict semifield \(\PosScalars\) is isomorphic to \(\PosReals\).
\end{proposition}
\begin{proof}
The proof of \cref{prop:posscalarssups} works \textit{mutatis mutandis}.
\end{proof}

\subsection{The scalars are the real or complex numbers}
Deducing that the field \(\Scalars\) is isomorphic to \(\Reals\) or \(\Comps\) is now purely a matter of algebra.

\begin{proposition}
\label{prop:reals}
There is an isomorphism of \(\SAScalars\) with \(\Reals\) that maps \(\PosScalars\) onto \(\PosReals\).
\end{proposition}

\begin{proof}
Let \(\varphi \colon \PosScalars \to \PosReals\) be the isomorphism in \cref{prop:pos_reals_concrete}. As \(\PosScalars\) contains all sums of squares of elements of \(\SAScalars\), we may define a map \(\psi \colon \SAScalars \to \Reals\) by
\[
    4\psi(a) = \varphi \paren[\big]{(a + 2)^2}  - \varphi\paren{a^2 + 4}
\]
for each \(a \in \SAScalars\). Then, for each \(a \in \PosScalars\), we have \(\psi(a) = \varphi(a)\) because
\[\varphi\paren[\big]{(a + 2)^2} = \varphi(a^2 + 4a + 4) = 4\varphi(a) + \varphi(a^2 + 4).\]
In particular, we have \(\psi(0) = 0\) and \(\psi(1) = 1\). Also, for all \(a, b \in \SAScalars\), we have
\begin{align*}
    4\psi(a + b)
    &= \varphi\paren[\big]{(a + b + 2)^2} - \varphi\paren[\big]{(a + b)^2 + 4},\\
    4\psi(a) + 4\psi(b)
    &= \varphi\paren[\big]{(a + 2)^2 + (b + 2)^2} - \varphi\paren[\big]{(a^2 + 4) +  (b^2 + 4)},\\\shortintertext{and}
    &\nqquad\nqquad(a + b + 2)^2 + (a^2 + 4) + (b^2 + 4)\\
    &= 2a^2 + 2b^2 + 2ab + 2a + 2b + 12\\
    &= (a+ 2)^2 + (b + 2)^2 + (a + b)^2 + 4,\\
\intertext{so \(\psi(a + b) = \psi(a) + \psi(b)\); we also have}
    16\psi(ab)
    &= \varphi\paren[\big]{4(ab + 2)^2} - \varphi\paren[\big]{4(a^2b^2 + 4)},\\
    16\psi(a)\psi(b)
    &= \varphi\paren[\big]{(a + 2)^2(b + 2)^2 + (a^2 + 4)(b^2 + 4)}\\
    &\qquad\qquad {} - \varphi\paren[\big]{(a + 2)^2(b^2 + 4) + (a^2 + 4)(b + 2)^2},\\\shortintertext{and}
    &\nqquad\nqquad 4(ab + 2)^2 + (a + 2)^2(b^2 + 4) + (a^2 + 4)(b + 2)^2\\
    &= 6a^2b^2 + 4a^2b + 4ab^2+ 8a^2 + 16ab + 8b^2 + 16a + 16b + 48\\
    &= (a + 2)^2(b + 2)^2 + (a^2 + 4)(b^2 + 4) + 4(a^2b^2 + 4),
\end{align*}
so \(\psi(ab) = \psi(a)\psi(b)\). Hence \(\psi\) is a ring homomorphism \(\SAScalars \to \Reals\) that extends \(\varphi\).

The map \(\upsilon \colon \Reals \to \SAScalars\) defined, for each \(a \in \Reals\), by
\[ 4\upsilon(a) = \varphi^{-1} \paren[\big]{(a + 2)^2}  - \varphi^{-1}\paren{a^2 + 4}\]
is similarly a ring homomorphism that extends \(\varphi^{-1}\). For each \(a \in \Reals\), we have
\begin{multline*}
    4\psi\upsilon(a) = \psi\paren[\big]{\varphi^{-1} \paren[\big]{(a + 2)^2} - \varphi^{-1} \paren{a^2 + 4}}
    = \psi\varphi^{-1} \paren[\big]{(a + 2)^2} - \psi\varphi^{-1} \paren{a^2 + 4}\\
    = \varphi\varphi^{-1} \paren[\big]{(a + 2)^2} - \varphi\varphi^{-1} \paren{a^2 + 4}
    = (a + 2)^2 - (a^2 + 4)
    = 4a,
\end{multline*}
and so \(\psi \upsilon = 1\). Similarly, we have \(\upsilon \psi = 1\). Hence \(\upsilon\) and \(\psi\) are mutually inverse.
\end{proof}

\begin{corollary}
There is an isomorphism of \(\Scalars\) with \(\Reals\) or \(\Comps\) that maps \(\PosScalars\) onto \(\PosReals\).
\end{corollary}

\begin{proof}
Combine \cref{prop:reals,lem:invol_reals_comps}.
\end{proof}

\subsection{Completeness} It remains unclear whether \cref{thm:main} still holds if \cref{axiom:colimits:mono} is replaced by \cref{axiom:colimits:epi}. The issue is establishing finite-dimensionality of the inner-product spaces corresponding to each object. This final subsection contains work towards proving \cref{prop:finite-dagger-biproduct} under our new assumptions.

For all objects \(A\), let \(\DaggerMonoCommaCat{\D}{A}\) be the full subcategory of \(\CommaCat{\D}{A}\) spanned by the objects \((X,x)\) where
% \(x \colon X \mathbin{\barrightarrowtail} A\)
\(x \colon \daggerMono{X}{A}\)
is dagger monic, and let \(\DaggerEpiCommaCat{A}{\D}\) be the full subcategory of \(\CommaCat{A}{\D}\) spanned by the objects \((X,x)\) where \(x \colon \daggerEpi{A}{X}\) is dagger epic. It follows from cancellativity that the morphisms in $\DaggerMonoCommaCat{\D}{A}$ and $\DaggerEpiCommaCat{A}{\D}$ are themselves dagger monic and dagger epic, respectively.

\begin{proposition}\label{prop:daggermonocolimit}
  For each object \(A\), the category \(\DaggerMonoCommaCat{\D}{A}\) has sequential colimits.
\end{proposition}
\begin{proof}
  Consider the following adjunction.
  \[
  \begin{tikzcd}[column sep=large, cramped]
  \CommaCat{\D}{A}
      \arrow[r, "\Coker", shift left=2]
      \arrow[r, left adjoint]
      \arrow[from=r, "\Ker", shift left=2]
      \&
  \CommaCat{A}{\D}
  \end{tikzcd}
  \]
  The functor \(\Ker\) maps each object \((X,x)\) of \(\CommaCat{A}{\D}\) to a chosen kernel of \(x\), and its action on morphisms is uniquely determined by universality of kernels. The functor \(\Coker\) is defined similarly. The unit \(\eta_{(X,x)} \colon (X,x) \to \Ker \Coker (X,x)\) of the adjunction is also uniquely determined by universality of kernels.

  As every normal monomorphism is a kernel of its cokernel, the adjunction is in fact idempotent, so~\cite[Theorem~3.8.7]{grandis:categorytheory} it factors as
  \[
  \begin{tikzcd}[column sep=large, cramped]
  \CommaCat{\D}{A}
      \arrow[r, "\Im", shift left=2]
      \arrow[r, left adjoint]
      \arrow[from=r, hook, shift left=2]
      \&
  \DaggerMonoCommaCat{\D}{A}
      \arrow[r, "\Coker", shift left=2]
      \arrow[r, left adjoint]
      \arrow[from=r, "\Ker", shift left=2]
      \&
  \DaggerEpiCommaCat{A}{\D}
      \arrow[r, hook, shift left=2]
      \arrow[r, left adjoint]
      \arrow[from=r, "\Coim", shift left=2]
      \&
  \CommaCat{A}{\D}
  \end{tikzcd}
  \]
  where the left, middle and right adjunctions are, respectively, a reflective subcategory inclusion, an equivalence of categories, and a coreflective subcategory inclusion.

  The canonical functor \(\CommaCat{A}{\D} \to \D\) creates connected colimits~\cite[Proposition~3.3.8]{riehl:categorycontext}, so in particular creates epimorphisms and sequential colimits. As \(\D\) has colimits of sequential diagrams of epimorphisms, so does \(\CommaCat{A}{\D}\). As \(\DaggerEpiCommaCat{A}{\D}\) is a reflective subcategory of \(\CommaCat{A}{\D}\), it has all colimits that \(\CommaCat{A}{\D}\) admits, formed by applying the reflector to the colimit~\cite[Proposition~4.5.15]{riehl:categorycontext}. In particular, it has sequential colimits. The result follows because the category \(\DaggerMonoCommaCat{\D}{A}\) is equivalent to \(\DaggerEpiCommaCat{A}{\D}\).
\end{proof}

\end{document}